\documentclass[10pt,reqno]{amsart}
\usepackage{amsmath, amssymb, amsthm}

\newcommand{\df}{\dfrac}

 \renewcommand{\a}{\alpha}
\renewcommand{\b}{\beta}
\newcommand{\e}{\epsilon}
\renewcommand{\d}{{\delta}}
\newcommand{\g}{\gamma}

\renewcommand{\(}{\left\(}
\renewcommand{\)}{\right\)}
\renewcommand{\[}{\left\[}
\renewcommand{\]}{\right\]}

\numberwithin{equation}{section}
 \theoremstyle{plain}
\newtheorem{theorem}{Theorem}[section]
\newtheorem{lemma}[theorem]{Lemma}

\newtheorem{corollary}[theorem]{Corollary}

   \makeatletter
\def\proof{\@ifnextchar[{\@oproof}{\@nproof}}
\def\@oproof[#1][#2]{\trivlist\item[\hskip\labelsep\textit{#2 Proof of\
#1.}~]\ignorespaces}
\def\@nproof{\trivlist\item[\hskip\labelsep\textit{Proof.}~]\ignorespaces}
\def\endproof{\qed\endtrivlist}
\makeatother

\begin{document}
\title[]{Overpartitions related to the mock theta function $\omega(q)$}

\author{George E. Andrews}
\address{Department of Mathematics, The Pennsylvania State University, University Park, PA 16802, USA} \email{gea1@psu.edu}

\author{Atul Dixit}
\address{Department of Mathematics, Indian Institute of Technology Gandhinagar, Near village Palaj, Gandhinagar, Gujarat 382355, India}\email{adixit@iitgn.ac.in}

\author{Daniel Schultz}
\address{Department of Mathematics, The Pennsylvania State University, University Park, PA 16802, USA}
\email{dps23@psu.edu}

\author{Ae Ja Yee}
\address{Department of Mathematics, The Pennsylvania State University, University Park, PA 16802, USA}
 \email{auy2@psu.edu}
 

\footnotetext[2]{2000 AMS Classification Numbers: Primary, 11P81; Secondary, 05A17}


\maketitle
\begin{abstract}
It was recently shown that $q\omega(q)$, where $\omega(q)$ is one of the third order mock theta functions, is the generating function of
$p_{\omega}(n)$, the number of partitions of a positive integer $n$ such that all odd parts are less than twice the smallest part. In this paper, we study the overpartition analogue of $p_{\omega}(n)$, and express its generating function in terms of a ${}_3\phi_{2}$ basic hypergeometric series and an infinite series involving little $q$-Jacobi polynomials. This is accomplished by obtaining a new seven parameter $q$-series identity  which generalizes a deep identity due to the first author as well as its generalization by R.P.~Agarwal. We also
derive two interesting congruences satisfied by the overpartition analogue, and some congruences satisfied by the associated smallest parts function.
\end{abstract}

\section{Introduction}\label{intro}
Since Ramanujan introduced mock theta functions in his last letter to Hardy in 1920, they have been the subject of intense study for many decades. Along with his third order mock theta function $f(q)$, there are many studies on the mock theta function
\begin{align*}
\omega(q):=\sum_{n=1}^{\infty} \frac{q^{2n^2+2n}}{(q;q^2)_{n+1}^2}
\end{align*}
in the literature \cite{bcgkmw}, \cite{jbko}, \cite{garthwaite}, \cite{waldherr}. Throughout the paper, we adopt the following $q$-series notation:
\begin{align*}
(a;q)_{n}&:=(1-a)(1-aq)\cdots (1-aq^{n-1}),\\
(a;q)_{\infty} &:=\lim_{n\to \infty} (a;q)_n, |q|<1.
\end{align*}
When the base is $q$, we sometimes use the short-hand notation $(a)_n:=(a;q)_n$, $(a)_{\infty}:=(a;q)_{\infty}$.

In a recent paper \cite{ady1}, the first, second and the fourth author discovered a new partition theoretic interpretation of $\omega(q)$, namely, the coefficient of $q^n$ in $q\,\omega(q)$ counts $p_{\omega}(n)$, the number of partitions of $n$ in which all odd parts are less than twice the smallest part, that is,
\begin{align*}
\sum_{n=1}^{\infty} p_{\omega}(n) q^n=\sum_{n=1}^{\infty} \frac{q^n}{(1-q^n)(q^{n+1};q)_{n} (q^{2n+2};q^2)_{\infty}}=q\omega(q).
\end{align*}
In the same paper they also studied the associated smallest parts function $\textup{spt}_{\omega}(n)$ whose generating function is given by
\begin{align*}
\sum_{n=1}^{\infty} \textup{spt}_{\omega}(n) q^n=\sum_{n=1}^{\infty} \frac{q^n}{(1-q^n)^2(q^{n+1};q)_{n} (q^{2n+2};q^2)_{\infty}}.
\end{align*}

In this paper we study the overpartition analogue of $p_{\omega}(n)$ and its associated smallest parts function. Overpartitions are ordinary partitions extended by allowing a possible overline designation on the first (or equivalently, the final) occurrence of a part. For instance, there are $8$ overpartitions of $3$, i.e., $3, \overline{3}, 2+1, 2+\overline{1}, \overline{2}+1, \overline{2}+\overline{1}, 1+1+1$, and $\overline{1}+1+1$. Throughout this paper, however, we consider overpartitions in which the smallest part is always overlined, and denote by $\overline{p}(n)$ the number of such overpartitions. For instance,  $\overline{p}(3)=4$ since there are $4$ such overpartitions of $3$, i.e., $\overline{3}, 2+\overline{1}, \overline{2}+\overline{1}$, and $\overline{1}+1+1$.  In an overpartition, a smallest part may or may not be overlined, so the number of overpartitions of $n$ is exactly twice of $\overline{p}(n)$. 

Since its introduction in \cite{corlove}, the overpartition construct has been very popular, and has led to a number of studies in $q$-series, partition theory, modular and mock modular forms.

As remarked earlier, in this paper we study the overpartition analogue of $p_{\omega}(n)$, namely $\overline{p}_{\omega}(n)$, which enumerates the number of overpartitions of $n$ such that all odd parts are less than twice the smallest part, and in which the smallest part is always overlined. It is clear that its generating function is given by
\begin{equation}\label{gfbpo}
\sum_{n=1}^{\infty}\overline{p}_{\omega}(n)q^n=\sum_{n=1}^{\infty}\frac{q^n(-q^{n+1};q)_n(-q^{2n+2};q^2)_{\infty}}{(1-q^n)(q^{n+1};q)_n(q^{2n+2};q^2)_{\infty}}.
\end{equation}
The series in \eqref{gfbpo} can be simplified to
\begin{align}
\sum_{n=1}^{\infty}\overline{p}_{\omega}(n)q^{n}
&=\frac{q(-q^2;q^2)_{\infty}}{(1-q)(q^2;q^2)_{\infty}}\sum_{n=0}^{\infty}\frac{(-q^3;q^2)_n(q;q)_n}{(q^3;q^2)_n(-q^2;q)_n}q^{n}\label{repre}\\
&=\frac{q(-q^2;q^2)_{\infty}}{(1-q)(q^2;q^2)_{\infty}}{}_{4}\phi_{3}\left(\begin{matrix}q,& q,&iq^{3/2},&-iq^{3/2}\\
&-q^2,& q^{3/2},& -q^{3/2}
\end{matrix}\, ;q, q \right),\nonumber
\end{align}
where the basic hypergeometric series ${}_{r+1}\phi_{r}$ is defined as
\begin{equation*}\label{bhs}
{}_{r+1}\phi_{r}\left(\begin{matrix} a_1, a_2, \ldots, a_{r+1}\\
  b_1,  b_2, \ldots, b_{r} \end{matrix}\,; q,
z \right) :=\sum_{n=0}^{\infty} \frac{(a_1;q)_n (a_2;q)_n \cdots (a_{r+1};q)_n}{(q;q)_n (b_1;q)_n \cdots (b_{r};q)_n} z^n.
\end{equation*}
Thus the generating function is essentially a nonterminating ${}_4\phi_3$. 
The problem of relating this generating function with familiar objects in the theory of basic hypergeometric series, modular forms and mock modular forms is quite difficult. In fact in order to transform it, we derive a new multi-parameter $q$-series identity, which generalizes a deep identity due to the first author, and its extension due to R.P.~Agarwal (see Theorem \ref{ext}). Basically we need its following variant for our case. 
\begin{theorem}\label{extvar}
Let the Gaussian polynomial be defined by
\begin{equation}\label{gaussian}
\left[\begin{matrix} n\\m\end{matrix}\right]=\left[\begin{matrix} n\\m\end{matrix}\right]_{q}:=
\begin{cases}
\displaystyle\frac{(q;q)_{n}}{(q;q)_m(q;q)_{n-m}},\hspace{2mm}\text{if}\hspace{1.5mm}0\leq m\leq n,\\
0,\hspace{2mm}\text{otherwise}.
\end{cases}
\end{equation}
Then, provided $\b,\d, f, t\neq q^{-j}, j\geq 0$, the following identity holds:
{\allowdisplaybreaks
\begin{align}\label{extandagar1}
&\sum_{n=0}^{\infty}\frac{(\a)_{n}(\g)_{n}(\e)_n}{(\b)_{n}(\d)_{n}(f)_{n}}t^n\nonumber\\
&=\frac{(\e, \g, \b/\a, q, \a t, q/(\a t), \d q/\b, fq/\b;q)_{\infty}}{(f, \d, q/\a, \b, \b/(\a t), \a tq/\b, \g q/\b, \e q/\b; q)_{\infty}}{}_{3}\phi_{2}\left(\begin{matrix} \frac{\a q}{\b},& \frac{\g q}{\b},& \frac{\e q}{\b} \\
& \frac{\d q}{\b},& \frac{fq}{\b} \end{matrix}\, ; q, t\right) \nonumber\\
&\quad+\left(1-\frac{q}{\b}\right)\frac{(\e, \g, t, \d q/\b, fq/\b;q)_{\infty}}{(f, \d, \a t/\b, \g q/\b, \e q/\b;q)_{\infty}}{}_{3}\phi_{2}\bigg(\begin{matrix}\frac{\a q}{\b},& \frac{\g q}{\b},&\frac{\e q}{\b}\\
& \frac{\d q}{\b},& \frac{fq}{\b}\end{matrix}\, ;q, t\bigg)
\bigg({}_{2}\phi_{1}\bigg(\begin{matrix}q,& \frac{q}{t}\\
& \frac{\b q}{\a t}\end{matrix}\, ;q, \frac{q}{\a} \bigg)-1\bigg)\nonumber\\
&\quad+\frac{(\e,\g;q)_{\infty}}{(f,\d;q)_{\infty}}\left(1-\frac{q}{\b}\right)\nonumber\\
&\quad\quad\times\sum_{n=0}^{\infty}\frac{(t)_n}{(q)_n(\a t/\b)_{n+1}}\left(\frac{q}{\b}\right)^n\sum_{p=0}^{n}\frac{(\a t/\b)_p}{(t)_p}\left(\frac{\b}{q}\right)^p\sum_{m=0}^{n}\left[\begin{matrix} n\\m\end{matrix}\right]\left(\frac{f}{\e}\right)_m\e^m\left(\frac{\d}{\g}\right)_{n-m}\g^{n-m},
\end{align}}
where we use the notation
\begin{equation*}
(a_1, a_2, \cdots, a_m;q)_n:=(a_1;q)_n(a_2;q)_n\cdots(a_m;q)_n.
\end{equation*}
\end{theorem}
This result is then specialized to obtain the following theorem which expresses the generating function in terms of a ${}_3\phi_{2}$ basic hypergeometric series and an infinite series involving the little $q$-Jacobi polynomial defined by \cite[Equation (3.1)]{aa}
\begin{equation}\label{little}
p_n(x;\a,\b:q)={}_{2}\phi_{1}\bigg(\begin{matrix}q^{-n},& \a\b q^{n+1}\\
&\a q\end{matrix}\, ;qx\bigg).
\end{equation}
\begin{theorem}\label{mt1}
The following identity holds for $|q|<1$:
\begin{align}\label{mainn}
\sum_{n=0}^{\infty}\frac{q^{n}(-q^3;q^2)_n(q;q)_n}{(q^3;q^2)_n(-q^2;q)_n}&=\frac{1}{2}\left(1-\frac{1}{q}\right)\frac{(q;q)_{\infty}^{2}}{(-q;q)_{\infty}^{2}}\sum_{n=0}^{\infty}\frac{(-1;q)_n(-q;q^2)_n}{(q;q^2)_n(q;q)_n}q^n\nonumber\\
&\quad+\frac{1}{q}\frac{(-q;q^2)_{\infty}}{(q^3;q^2)_{\infty}}\sum_{n=0}^{\infty}\frac{(q;q^2)_n(-q)^n}{(-q;q^2)_{n}(1+q^{2n})}p_{2n}(-1;q^{-2n-1},-1:q).
\end{align}
Hence the generating function of $\overline{p}_{\omega}(n)$ is given by
\begin{align}\label{fgf}
\sum_{n=1}^{\infty}\overline{p}_{\omega}(n)q^n&=-\frac{1}{2}\frac{(q;q)_{\infty}(q;q^2)_{\infty}}{(-q;q)_{\infty}(-q;q^2)_{\infty}}\sum_{n=0}^{\infty}\frac{(-1;q)_n(-q;q^2)_n}{(q;q^2)_n(q;q)_n}q^n\nonumber\\
&\quad+\frac{(-q;q)_{\infty}}{(q;q)_{\infty}}\sum_{n=0}^{\infty}\frac{(q;q^2)_n(-q)^n}{(-q;q^2)_{n}(1+q^{2n})}p_{2n}(-1;q^{-2n-1},-1:q).
\end{align}
\end{theorem}
The series involving the little $q$-Jacobi polynomials on the right side of \eqref{mainn} satisfies a nice congruence modulo $4$ given below.
\begin{theorem}\label{congmod4}
The following congruence holds:
\begin{equation}\label{mod4}
\frac{1}{q}\frac{(-q;q^2)_{\infty}}{(q^3;q^2)_{\infty}}\sum_{n=0}^{\infty}\frac{(q;q^2)_n(-q)^n}{(-q;q^2)_{n}(1+q^{2n})}p_{2n}(-1;q^{-2n-1},-1:q)\equiv\frac{1}{2q}-\frac{1}{2}\pmod{4}.
\end{equation}
\end{theorem}
\noindent
This suggests that this series might be further linked to some important objects in the literature. As of now this has remained elusive for us though.

The overpartition function $\overline{p}_{\omega}(n)$ satisfies some nice congruences. Indeed the two congruences in the following theorem will be proved in Section \ref{pwcong}.
\begin{theorem}\label{pw4}
We have
\begin{align}
\overline{p}_{\omega}(4n+3)\equiv0\pmod{4},\label{pw4cong1}\\
\overline{p}_{\omega}(8n+6)\equiv0\pmod{4}.\label{pw4cong2}
\end{align}
\end{theorem}

\noindent
Bringmann, Lovejoy, and Osburn \cite{bringmannlovejoyosburn1, bringmannlovejoyosburn2} defined $\overline{\textup{\textup{spt}}}(n)$ as the number of smallest parts in the overpartitions of $n$ and showed that $\overline{\textup{\textup{spt}}}(n)$ is a quasimock theta function (see \cite[p.~3--4]{bringmannlovejoyosburn2} for the definition) satisfying simple Ramanujan-type congruences, for instance,
\begin{align*}
\overline{\textup{spt}}(3n)\equiv 0 \pmod{3}. 
\end{align*}
In this paper, we study $\overline{\textup{\textup{spt}}}_{\omega}(n)$, the number of smallest parts in the overpartitions of $n$ in which the smallest part is always overlined and all odd parts are less than twice the smallest part.  
By its definition we see that the generating function of  $\overline{\textup{\textup{spt}}}_{\omega}(n)$ is given by
 \begin{align} \label{sptbaromega}
 \sum_{n=1}^{\infty} \overline{\textup{\textup{spt}}}_{\omega}(n) q^n=\sum_{n=1}^{\infty} \frac{q^{n}(-q^{n+1};q)_{n} (-q^{2n+2};q^2)_{\infty}}{(1-q^n)^2 (q^{n+1};q)_n (q^{2n+2};q^2)_{\infty}}.
 \end{align}
The smallest parts function $\overline{\textup{\textup{spt}}}_{\omega}(n)$ seems to carry arithmetic properties analogous to those of $\overline{\textup{\textup{spt}}}_2(n)$, where 
$\overline{\textup{\textup{spt}}}_2(n)$ counts the number of smallest parts in the overpartitions of $n$ with smallest parts even. It is known \cite{bringmannlovejoyosburn1} that 
\begin{align}
  \overline{\textup{\textup{spt}}}_2(3n) &\equiv 0 \pmod{3}, \label{eq3n0}\\
  \overline{\textup{\textup{spt}}}_{2}(3n+1) &\equiv 0 \pmod{3}, \label{eq3n1}\\
  \overline{\textup{\textup{spt}}}_{2}(5n+3)&\equiv 0 \pmod{5}. \label{eq5n3}
 \end{align}
The following are the main congruences satisfied by $\overline{\textup{\textup{spt}}}_{\omega}(n)$:
 
\begin{theorem} \label{thm1}
We have 
 \begin{align}
  \overline{\textup{\textup{spt}}}_\omega(3n) &\equiv 0 \pmod{3}\text{,} \label{equ_3n0} \\
  \overline{\textup{\textup{spt}}}_{\omega}(3n+2) &\equiv 0 \pmod{3}\text{,} \label{equ_3n2}\\
  \overline{\textup{\textup{spt}}}_{\omega}(10n+6)&\equiv 0 \pmod{5}\text{,} \label{equ_10n6}\\
	\overline{\textup{\textup{spt}}}_{\omega}(6n+5) &\equiv 0 \pmod{6}. \label{eq6n5}
 \end{align}
\end{theorem}

There are further congruences that both $\overline{\textup{\textup{spt}}}(n)$ and $\overline{\textup{\textup{spt}}}_{\omega}(n)$ satisfy:
\begin{theorem} \label{thm2}
For any positive integer $n$,
\begin{align}
\overline{\textup{\textup{spt}}}_{\omega}(n) \equiv \overline{\textup{\textup{spt}}}(n)&\equiv \begin{cases} 1 \pmod{2} & \text{ if $n= k^2$ or $ 2 k^2$ for some $k$,} \\
0 \pmod{2} & \text{otherwise}. \end{cases} \label{equ_nm2}
\end{align}
\end{theorem}

\begin{theorem} \label{thm3}
For any positive integer $n$,
\begin{align}
\overline{\textup{\textup{spt}}}(7n)&\equiv \overline{\textup{\textup{spt}}}(\tfrac{n}{7})  \pmod{4}\text{,} \label{equ_7nm4} \\
\overline{\textup{\textup{spt}}}_{\omega}(7n)&\equiv \overline{\textup{\textup{spt}}}_{\omega}(\tfrac{n}{7})  \pmod{4}\text{,} \label{equ_w7nm4}
\end{align}
where we follow the convention that $\overline{\textup{\textup{spt}}}(x)=\overline{\textup{\textup{spt}}}_{\omega}(x)=0$ if $x$ is not a positive integer.
\end{theorem}

This paper is organized as follows. In Section~\ref{prelim}, we recall some basic facts and theorems that are used in the sequel. Section \ref{bpw} is devoted to finding an alternate representation for the generating function of  $\overline{p}_{\omega}(n)$ in terms of a ${}_3\phi_{2}$ basic hypergeometric series and an infinite series involving the little $q$-Jacobi polynomials. A congruence modulo $4$, satisfied by the latter series, is also obtained in this section. In Section \ref{pwcong}, we give a proof of the congruences modulo $4$ satisfied by $\overline{p}_{\omega}(n)$. We recall some facts about $\overline{\textup{spt}}(n)$ and $\overline{\textup{spt}}_2(n)$ in Section \ref{sptwdf} and represent the generating function of $\overline{\textup{spt}}_{\omega}(n)$ in terms of those of these functions. In Section~\ref{sptwcong}, we prove the congruences modulo $3, 5$ and $6$ given in Theorem~\ref{thm1} based on these representations. Lastly we prove Theorems \ref{thm2} and \ref{thm3} in Section \ref{congsptsptw}. 

\section{Preliminaries}\label{prelim}
We collect below the important facts and theorems in the literature on $q$-series and partitions in this section. First of all, we assume throughout the paper that $|q|<1$. 
The most fundamental theorem in the literature is the $q$-binomial theorem given for $|z|<1$ by \cite[p.~17, Equation (2.2.1)]{gea1998}
\begin{equation}\label{qbin}
\sum_{n=0}^{\infty}\frac{(a;q)_{n}z^n}{(q;q)_n}=\frac{(az;q)_{\infty}}{(z;q)_{\infty}}.
\end{equation}
For $|z|<1$ and $|b|<1$, Heine's transformation \cite[p.~19, Corollary 2.3]{gea1998} is given by
\begin{equation}\label{heine}
{}_2\phi_{1}\left(\begin{matrix}a,& b\\
&c\end{matrix}\, ;q, z\right)=\frac{(b, az;q)_{\infty}}{(c, z;q)_{\infty}}{}_2\phi_{1}\left(\begin{matrix}c/b,& z\\
&az\end{matrix}\, ;q, b\right),
\end{equation}
and we also note Bailey's ${}_{10}\phi_{9}$ transformation \cite[Equation (2.10)]{andrews1984}, \cite[Equation (6.3)]{bailey}
{\allowdisplaybreaks\begin{align}
&\lim_{N\to\infty}{}_{10}\phi_{9}\left(\begin{matrix} a,& q^2\sqrt{a},& -q^2\sqrt{a},& p_1, &
p_1q, &p_2, & p_2q, & f, & q^{-2N}, & q^{-2N+1}\\
 &\sqrt{a}, & -\sqrt{a},& \df{a
q^2}{p_1}, & \df{a q}{p_1}, & \df{a q^2}{p_2}, & \df{a q}{p_2}, & \df{aq^2}{f}, & aq^{2N+2}, & aq^{2N+1}
\end{matrix}\,; q^2,
 \df{a^3q^{4N+3}}{p_1^{2}p_2^{2}f}\right) \notag \\
& \quad =\df{(a q;q)_{\infty}\left(\df{a q}{p_1p_2};q\right)_{\infty}}
{\left(\df{a q}{p_1};q\right)_{\infty}\left(\df{a q}{p_2};q\right)_{\infty}}\sum_{n=0}^{\infty}\frac{(p_1;q)_n(p_2;q)_n\left(\frac{aq}{f};q^2\right)_n}{(q;q)_n(aq;q^2)_n\left(\frac{aq}{f};q\right)_n}\left(\frac{aq}{p_1p_2}\right)^n. \label{baitra}
\end{align}}
Finally we note a transformation for ${}_2\psi_{2}$ due to Bailey \cite[p.~148, Exer. 5.11]{gasper}:
\begin{align}
{_2} \psi_2 \left(\begin{matrix} e, f\\ \frac{aq}{c}, \frac{aq}{d} \end{matrix}\, ; q, \frac{aq}{ef}\right)=\frac{(\frac{q}{c}, \frac{q}{d}, \frac{aq}{e}, \frac{aq}{f};q)_{\infty}}{(aq, \frac{q}{a}, \frac{aq}{cd}, \frac{aq}{ef};q)_{\infty}}\sum_{n=-\infty}^{\infty} \frac{(1-aq^{2n})(c,d,e,f;q)_n}{(1-a)(\frac{aq}{c}, \frac{aq}{d}, \frac{aq}{e}, \frac{aq}{f};q)_n} \left(\frac{qa^3}{cdef}\right)^n q^{n^2}, \label{9}
\end{align}
where ${}_r\psi_{r}$ is the basic bilateral hypergeometric series defined by \cite[p.~137, Equation (5.1.1)]{gasper}
\begin{equation*}
{}_r\psi_{r}\left(\begin{matrix} a_1, a_2, \ldots, a_{r}\\
  b_1,  b_2, \ldots, b_{r} \end{matrix}\,; q,
z \right):=\sum_{n=-\infty}^{\infty}\frac{(a_1, a_2, \cdots, a_r;q)_n}{(b_1, b_2, \cdots, b_r;q)_n}z^n.
\end{equation*}
\section{The generating function of $\overline{p}_{\omega}(n)$}\label{bpw}
First, we recall that $\overline{p}_{\omega}(n)$ counts the number of overpartitions of $n$ such that all odd parts are less than twice the smallest part, and in which the smallest part is always overlined. None of the already existing identities from the theory of basic hypergeometric series seems to be capable of handling its generating function. Hence we devise a new $q$-series identity consisting of seven parameters that transforms \eqref{repre} into a ${}_3\phi_{2}$ and an infinite series involving little $q$-Jacobi polynomials defined in \eqref{little}. The motivation and the need for devising such an identity is now given.

In the proof of the representation of the generating function of $p_{\omega}(n)$ in terms of the third order mock theta function $\omega(q)$ \cite[Theorem 3.1]{ady1}, the following four parameter $q$-series identity due to the first author \cite[Theorem 1]{gea90} played an instrumental role.
\begin{align}
\sum_{n=0}^{\infty} \frac{(B;q)_n (-Abq;q)_n q^n}{(-aq;q)_n (-bq;q)_n}&=\frac{-a^{-1} (B;q)_{\infty} (-Abq;q)_{\infty}}{(-bq;q)_{\infty} (-aq;q)_{\infty}} \sum_{m=0}^{\infty} \frac{(A^{-1};q)_m \left(\frac{Abq}{a}\right)^m}{\left(-\frac{B}{a};q\right)_{m+1}}\nonumber\\
&\quad+(1+b) \sum_{m=0}^{\infty} \frac{(-a^{-1};q)_{m+1} \left( -\frac{ABq}{a};q\right)_{m} (-b)^m}{\left(-\frac{B}{a};q\right)_{m+1} \left(\frac{Abq}{a};q\right)_{m+1}}. \label{gea90_thm1}
\end{align}

Agarwal \cite[Equation (3.1)]{agar1} obtained the following `mild' extension/generalization of \eqref{gea90_thm1} in the sense that we get \eqref{gea90_thm1} from the following identity when $t=q$.
\begin{align}\label{me}
&\sum_{n=0}^{\infty}\frac{(\a)_{n}(\g)_{n}}{(\b)_{n}(\d)_{n}}t^n\nonumber\\
&=\frac{(q/(\a t), \g, \a t, \b/\a, q;q)_{\infty}}{(\b/(\a t), \d, t, q/\a, \b;q)_{\infty}}{}_{2}\phi_{1}\bigg(\begin{matrix}\d/\g,& t\\
&q\a t/\b\end{matrix}\, ;q, \g q/\b\bigg)\nonumber\\
&\quad+\frac{(\g)_{\infty}}{(\d)_{\infty}}\left(1-\frac{q}{\b}\right)\sum_{m=0}^{\infty}\frac{(\d/\g)_{m}(t)_m}{(q)_m(\a t/\b)_{m+1}}(q\g/\b)^m\left({}_{2}\phi_{1}\bigg(\begin{matrix}q,& q/t\\
&q\b/(\a t)\end{matrix}\, ;q, q/\a\bigg)-1\right)\nonumber\\
&\quad+\frac{(\g)_{\infty}}{(\d)_{\infty}}\left(1-\frac{q}{\b}\right)\sum_{p=0}^{\infty}\frac{\g^p(\d/\g)_p}{(q)_p}\sum_{m=0}^{\infty}\frac{(\d q^p/\g)_m(tq^p)_m}{(q^{1+p})_m(\a tq^p/\b)_{m+1}}(q\g/\b)^m.
\end{align}
Since the right side of \eqref{repre} involves three $q$-shifted factorials (with base $q$) in the numerator as well as in the denominator of its summand, we need to first generalize \eqref{me}. Indeed, a generalization of \eqref{me} will be given below. However, we shall first prove its variant, namely Theorem \ref{extvar}, that we need for our purpose.
\begin{proof}[Theorem \textup{\ref{extvar}}][]
Let
\begin{align}\label{sabb}
S:=S(\a, \b, \g, \d, \e, f;q;t):=\sum_{n=0}^{\infty}\frac{(\a)_{n}(\g)_{n}(\e)_n}{(\b)_{n}(\d)_{n}(f)_{n}}t^n.
\end{align}
Then by an application of the $q$-binomial theorem \eqref{qbin},
\begin{align}\label{s}
S&=\frac{(\e)_{\infty}}{(f)_{\infty}}\sum_{n=0}^{\infty}\frac{(\a)_{n}(\g)_{n}(fq^{n})_{\infty}}{(\b)_{n}(\d)_{n}(\e q^n)_{\infty}}t^n\nonumber\\
&=\frac{(\e)_{\infty}}{(f)_{\infty}}\sum_{n=0}^{\infty}\frac{(\a)_{n}(\g)_{n}t^n}{(\b)_{n}(\d)_{n}}\sum_{m=0}^{\infty}\frac{(f/\e)_m}{(q)_m}(\e q^n)^m\nonumber\\
&=\frac{(\e)_{\infty}}{(f)_{\infty}}\sum_{m=0}^{\infty}\frac{(f/\e)_m\e^m}{(q)_m}\sum_{n=0}^{\infty}\frac{(\a)_{n}(\g)_{n}}{(\b)_{n}(\d)_{n}}(tq^m)^n\nonumber\\
&=\frac{(\e)_{\infty}}{(f)_{\infty}}\sum_{m=0}^{\infty}\frac{(f/\e)_m\e^m}{(q)_m}\bigg\{\frac{(q^{1-m}/(\a t), \g, \a tq^m, \b/\a, q;q)_{\infty}}{(\b q^{-m}/(\a t), \d, tq^m, q/\a, \b;q)_{\infty}}{}_{2}\phi_{1}\bigg(\begin{matrix}\d/\g,& tq^m\\
&\a tq^{m+1}/\b\end{matrix}\, ;q, \g q/\b\bigg)\nonumber\\
&\quad\quad\quad\quad\quad\quad\quad\quad\quad+\frac{(\g)_{\infty}}{(\d)_{\infty}}\frac{\left(1-q/\b\right)}{(1-\a tq^m/\b)}\sum_{k=0}^{\infty}\frac{(\d/\g)_k(\a tq^m/\b)_k\g^k}{(q)_k(\a tq^{m+1}/\b)_k}\sum_{r=0}^{\infty}\frac{(q^{1-k-m}/t)_r}{(\b q^{1-k-m}/(\a t))_r}\left(\frac{q}{\a}\right)^r\bigg\}\nonumber\\
&=:\frac{(\e)_{\infty}}{(f)_{\infty}}\left(\frac{(\g,\b/\a, q;q)_{\infty}}{(\d, q/\a, \b;q)_{\infty}}V_1+V_2\right),
\end{align}
where in the penultimate step, we used \eqref{me} in the form given in \cite[Equation (3.2)]{agar1}. Here
\begin{align}
V_1&:=\sum_{m=0}^{\infty}\frac{(f/\e)_m(q^{1-m}/(\a t),\a tq^m;q)_{\infty}\e^m}{(q)_m(\b q^{-m}/(\a t), tq^m;q)_{\infty}}{}_{2}\phi_{1}\bigg(\begin{matrix}\d/\g,& tq^m\\
&\a tq^{m+1}/\b\end{matrix}\, ;q, \g q/\b\bigg),\nonumber\\
V_2&:=\frac{(\g)_{\infty}}{(\d)_{\infty}}\left(1-\frac{q}{\b}\right)\sum_{m=0}^{\infty}\frac{\left(\frac{f}{\e}\right)_m\e^m}{(q)_m\left(1-\frac{\a tq^m}{\b}\right)}\sum_{k=0}^{\infty}\frac{\left(\frac{\d}{\g}\right)_k\left(\frac{\a tq^m}{\b}\right)_k\g^k}{(q)_k\left(\frac{\a tq^{m+1}}{\b}\right)_k}\sum_{r=0}^{\infty}\frac{\left(\frac{q^{1-k-m}}{t}\right)_r}{\left(\frac{\b q^{1-k-m}}{\a t}\right)_r}\left(\frac{q}{\a}\right)^r.\label{vv12}
\end{align}
Next, using Heine's transformation \eqref{heine} in the second step below, we see that
\begin{align}
V_1&=\frac{(\a t)_{\infty}(q/(a t))_{\infty}}{( t)_{\infty}(\b/(a t))_{\infty}}\sum_{m=0}^{\infty}\frac{(f/\e)_m(t)_m}{(\a tq/\b)_m(q)_m}\left(\frac{q\e}{\b}\right)^m{}_{2}\phi_{1}\bigg(\begin{matrix}\d/\g,& tq^m\\
&\a tq^{m+1}/\b\end{matrix}\, ;q, \g q/\b\bigg)\nonumber\\
&=\frac{(\a t)_{\infty}(q/(a t))_{\infty}(\d q/\b)_{\infty}}{( \a tq/\b)_{\infty}(\g q/\b)_{\infty}(\b/(a t))_{\infty}}\sum_{m=0}^{\infty}\frac{(f/\e)_m}{(q)_m}\left(\frac{q\e}{\b}\right)^m{}_{2}\phi_{1}\bigg(\begin{matrix}\a q/\b,& \g q/\b\\
&\d q/\b\end{matrix}\, ;q, tq^m\bigg)\nonumber\\
&=\frac{(\a t)_{\infty}(q/(a t))_{\infty}(\d q/\b)_{\infty}}{( \a tq/\b)_{\infty}(\g q/\b)_{\infty}(\b/(a t))_{\infty}}\sum_{k=0}^{\infty}\frac{(\a q/\b)_k(\g q/\b)_k(fq^{k+1}/\b)_{\infty}t^k}{(\d q/\b)_k(q)_k(\e q^{k+1}/\b)_{\infty}},
\end{align}
where in the last step we used \eqref{qbin} after interchanging the order of summation. Hence
\begin{align}\label{v1}
V_1=\frac{(\a t, q/(\a t), \d q/\b, fq/\b;q)_{\infty}}{(\b/(\a t), \a tq/\b, \g q/\b, \e q/\b;q)_{\infty}}{}_{3}\phi_{2}\bigg(\begin{matrix}\a q/\b,& \g q/\b, & \e q/\b\\
&\d q/\b, fq/\b\end{matrix}\, ;q, t\bigg).
\end{align}
Let us now consider $V_2$. Since
{\allowdisplaybreaks\begin{align}\label{v2int}
\sum_{r=0}^{\infty}\frac{(q^{1-k-m}/t)_r}{(\b q^{1-k-m}/(\a t))_r}\left(\frac{q}{\a}\right)^r&=\frac{(t)_{m+k}}{(\a t/\b)_{m+k}}\left(\frac{q}{\b}\right)^{m+k}\sum_{r=0}^{\infty}\frac{(\a t/\b)_{m+k-r}}{(t)_{m+k-r}}\left(\frac{\b}{q}\right)^{m+k-r}\nonumber\\
&=\frac{(t)_{m+k}}{\left(\frac{\a t}{\b}\right)_{m+k}}\left(\frac{q}{\b}\right)^{m+k}\left(\sum_{p=1}^{\infty}\frac{\left(\frac{q}{t}\right)_p}{\left(\frac{\b q}{\a t}\right)_p}\left(\frac{q}{\a}\right)^p+\sum_{p=0}^{m+k}\frac{(\a t/\b)_p}{(t)_p}\left(\frac{\b}{q}\right)^p\right),
\end{align}}
we find that
\begin{equation}\label{v21s}
V_2=\frac{(\g)_{\infty}}{(\d)_{\infty}}\left(1-\frac{q}{\b}\right)(V_3+V_3^{*}),
\end{equation}
where
\begin{align}\label{v345ms}
V_3&:=\sum_{m=0}^{\infty}\frac{(f/\e)_m(t)_m(\e q/\b)^m}{(q)_m(\a t/\b)_m(1-\a tq^m/\b)}\sum_{k=0}^{\infty}\frac{(\d/\g)_k(tq^m)_k(\g q/\b)^k}{(q)_k(\a tq^{m+1}/\b)_k}\sum_{p=1}^{\infty}\frac{(q/t)_p}{(\b q/(\a t))_p}\left(\frac{q}{\a}\right)^p,\nonumber\\
V_3^{*}&:=\sum_{m=0}^{\infty}\frac{(f/\e)_m(\e q/\b)^m}{(q)_m(1-\a tq^m/\b)}\sum_{k=0}^{\infty}\frac{(\d/\g)_k(t)_{m+k}(\a tq^m/\b)_k(\g q/\b)^k}{(q)_k(\a tq^{m+1}/\b)_k(\a t/\b)_{m+k}}\sum_{p=0}^{m+k}\frac{(\a t/\b)_p}{(t)_p}\left(\frac{\b}{q}\right)^p.
\end{align}
Consider $V_3$. Again using Heine's transformation \eqref{heine} for the middle series followed by \eqref{qbin}, we have 
\begin{align}\label{v3}
V_3&=\frac{(t, \d q/\b, fq/\b;q)_{\infty}}{(\a t/\b, \g q/\b, \e q/\b;q)_{\infty}}{}_{3}\phi_{2}\bigg(\begin{matrix}\a q/\b,& \g q/\b,&\e q/\b\\
&\d q/\b,& fq/\b\end{matrix}\, ;q, t\bigg)\bigg({}_{2}\phi_{1}\bigg(\begin{matrix}q,& q/t\\
&\b q/(\a t)\end{matrix}\, ;q, q/\a\bigg)-1\bigg).
\end{align}
Also,
\begin{align}\label{v3s}
V_3^{*}=\sum_{n=0}^{\infty}\frac{(t)_n}{(q)_n(\a t/\b)_{n+1}}\left(\frac{q}{\b}\right)^n\sum_{p=0}^{n}\frac{(\a t/\b)_p}{(t)_p}\left(\frac{\b}{q}\right)^p\sum_{m=0}^{n}\left[\begin{matrix} n\\m\end{matrix}\right]\left(\frac{f}{\e}\right)_m\e^m\left(\frac{\d}{\g}\right)_{n-m}(\g)^{n-m}.
\end{align}
Finally from \eqref{s}, \eqref{v1}, \eqref{v21s}, \eqref{v3} and \eqref{v3s}, we arrive at \eqref{extandagar1}.
\end{proof}
We now give the aforementioned generalization of \eqref{me} which can also be viewed as a corollary of Theorem \ref{extvar}.
\begin{theorem}\label{ext}
Provided $\b,\d, f, t\neq q^{-j}, j\geq 0$, the following identity holds:
{\allowdisplaybreaks\begin{align}\label{extandagar}
&\sum_{n=0}^{\infty}\frac{(\a)_{n}(\g)_{n}(\e)_n}{(\b)_{n}(\d)_{n}(f)_{n}}t^n\nonumber\\
&=\frac{(\e, \g, \b/\a, q, \a t, q/(\a t), \d q/\b, fq/\b;q)_{\infty}}{(f, \d, q/\a, \b, \b/(\a t), \a tq/\b, \g q/\b, \e q/\b; q)_{\infty}}{}_{3}\phi_{2}\left(\begin{matrix}\a q/\b,& \g q/\b,&\e q/\b\\
&\d q/\b,& fq/\b\end{matrix}\, ;q, t\right)\nonumber\\
&\quad+\left(1-\frac{q}{\b}\right)\frac{(\e, \g, t, \d q/\b, fq/\b;q)_{\infty}}{(f, \d, \a t/\b, \g q/\b, \e q/\b;q)_{\infty}}{}_{3}\phi_{2}\bigg(\begin{matrix}\a q/\b,& \g q/\b,&\e q/\b\\
&\d q/\b,& fq/\b\end{matrix}\, ;q, t\bigg)\nonumber\\
&\quad\quad\times\bigg({}_{2}\phi_{1}\bigg(\begin{matrix}q,& q/t\\
&\b q/(\a t)\end{matrix}\, ;q/\a\bigg)-1\bigg)\nonumber\\
&\quad+\left(1-\frac{q}{\b}\right)\frac{(\e, \g, t, fq/\b;q)_{\infty}}{(f, \d, \a t/\b, \e q/\b;q)_{\infty}}\nonumber\\
&\quad\quad\times\sum_{p=0}^{\infty}\frac{(\d/\g)_p(\a t/\b)_p\g^p}{(t)_p(q)_p}\sum_{k=0}^{\infty}\frac{(\d q^p/\g)_k(q\g/\b)^k}{(q^{1+p})_k}{}_{2}\phi_{1}\bigg(\begin{matrix}\a q/b,& \e q/\b\\
&fq/\b\end{matrix}\, ;q, tq^{k+p}\bigg)\nonumber\\
&\quad+\left(1-\frac{q}{\b}\right)\frac{(\e, \g;q)_{\infty}}{(f, \d;q)_{\infty}}\sum_{p=1}^{\infty}\frac{(f/\e)_p\e^p}{(q)_p}\sum_{k=0}^{\infty}\frac{(\d/\g)_k\g^k}{(q)_k}\sum_{m=0}^{\infty}\frac{(fq^p/\e)_m(tq^{p+k})_m}{(q^{1+p})_m(\a tq^{p+k}/\b)_{m+1}}(\e q/\b)^{m}.
\end{align}}
\end{theorem}
\begin{proof}
Write $V_3^{*}$ in \eqref{v345ms} as
\begin{equation}\label{v3si}
V_3^{*}=V_4+V_5,
\end{equation}
where
\begin{align}\label{v345m}
V_4&=\sum_{m=0}^{\infty}\frac{(f/\e)_m(\e q/\b)^m}{(q)_m(1-\a tq^m/\b)}\sum_{k=0}^{\infty}\frac{(\d/\g)_k(t)_{m+k}(\a tq^m/\b)_k(\g q/\b)^k}{(q)_k(\a tq^{m+1}/\b)_k(\a t/\b)_{m+k}}\sum_{p=0}^{k}\frac{(\a t/\b)_p}{(t)_p}\left(\frac{\b}{q}\right)^p,\nonumber\\
V_5&=\sum_{m=0}^{\infty}\frac{(f/\e)_m(\e q/\b)^m}{(q)_m(1-\a tq^m/\b)}\sum_{k=0}^{\infty}\frac{(\d/\g)_k(t)_{m+k}(\a tq^m/\b)_k(\g q/\b)^k}{(q)_k(\a tq^{m+1}/\b)_k(\a t/\b)_{m+k}}\sum_{p=1}^{m}\frac{(\a t/\b)_{k+p}}{(t)_{k+p}}\left(\frac{\b}{q}\right)^{k+p}.
\end{align}
Note that $V_4$ can be written as
{\allowdisplaybreaks\begin{align}\label{v4}
V_4&=\sum_{m=0}^{\infty}\frac{(f/\e)_m(t)_m(\e q/\b)^m}{(q)_m(\a t/\b)_{m+1}}\sum_{p=0}^{\infty}\frac{(\a t/\b)_p}{(t)_p}\left(\frac{\b}{q}\right)^p\sum_{k=p}^{\infty}\frac{(\d/\g)_k(tq^m)_k}{(q)_k(\a tq^{m+1}/\b)_k}\left(\frac{\g q}{\b}\right)^k\nonumber\\
&=\sum_{m=0}^{\infty}\frac{(f/\e)_m(t)_m(\e q/\b)^m}{(q)_m(\a t/\b)_{m+1}}\sum_{p=0}^{\infty}\frac{(\a t/\b)_p(\d/\g)_p(tq^m)_p\g^p}{(t)_p(\a tq^{m+1}/\b)_p(q)_p}\sum_{k=0}^{\infty}\frac{(\d q^p/\g)_k(tq^{m+p})_k(\g q/\b)^k}{(q^{p+1})_k(\a tq^{m+p+1}/\b)_k}\nonumber\\
&=\sum_{p=0}^{\infty}\frac{(\d/\g)_p\g^p}{(q)_p}\sum_{m=0}^{\infty}\frac{(f/\e)_m(tq^p)_m}{(q)_m(\a tq^p/\b)_{m+1}}\left(\frac{\e q}{\b}\right)^m\sum_{k=0}^{\infty}\frac{(\d q^p/\g)_k(tq^{m+p})_k(\g q/\b)^k}{(q^{p+1})_k(\a tq^{m+p+1}/\b)_k}\nonumber\\
&=\sum_{p=0}^{\infty}\frac{(\d/\g)_p\g^p}{(q)_p}\sum_{k=0}^{\infty}\frac{(\d q^p/\g)_k(tq^p)_k}{(q^{p+1})_k(\a tq^p/\b)_{k+1}}\left(\frac{\g q}{\b}\right)^k\sum_{m=0}^{\infty}\frac{(f/\e)_m(tq^{p+k})_m}{(q)_m(\a tq^{p+k+1}/\b)_m}\left(\frac{\e q}{\b}\right)^m\nonumber\\
&=\frac{(fq/\b)_{\infty}(t)_{\infty}}{(\a t/\b)_{\infty}(\e q/\b)_{\infty}}\sum_{p=0}^{\infty}\frac{(\d/\g)_p(\a t/\b)_p\g^p}{(t)_p(q)_p}\sum_{k=0}^{\infty}\frac{(\d q^p/\g)_k(q\g/\b)^k}{(q^{1+p})_k}{}_{2}\phi_{1}\bigg(\begin{matrix}\a q/b,& \e q/\b\\
&fq/\b\end{matrix}\, ;q, tq^{k+p}\bigg),
\end{align}}%
where in the last step, we used \eqref{heine} to transform the innermost series. Lastly, $V_5$ can be simplified to
\begin{align}\label{v5}
V_5&=\sum_{k=0}^{\infty}\frac{(\d/\g)_k\g^k}{(q)_k}\sum_{m=1}^{\infty}\frac{(f/\e)_m(tq^k)_m(\e q/\b)^m}{(q)_m(\a tq^k/\b)_{m+1}}\sum_{p=1}^{m}\frac{(\a tq^k/\b)_{p}}{(tq^k)_p}\left(\frac{\b}{q}\right)^p\nonumber\\
&=\sum_{k=0}^{\infty}\frac{(\d/\g)_k\g^k}{(q)_k}\sum_{p=1}^{\infty}\frac{(\a tq^k/\b)_{p}}{(tq^k)_p}\left(\frac{\b}{q}\right)^p\sum_{m=0}^{\infty}\frac{(f/\e)_{m+p}(tq^k)_{m+p}}{(q)_{m+p}(\a tq^k/\b)_{m+p+1}}\left(\frac{\e q}{\b}\right)^{m+p}\nonumber\\
&=\sum_{p=1}^{\infty}\frac{(f/\e)_p\e^p}{(q)_p}\sum_{k=0}^{\infty}\frac{(\d/\g)_k\g^k}{(q)_k}\sum_{m=0}^{\infty}\frac{(fq^p/\e)_m(tq^{p+k})_m}{(q^{1+p})_m(\a tq^{p+k}/\b)_{m+1}}\left(\frac{\e q}{\b}\right)^{m}.
\end{align}
Now \eqref{s}, \eqref{v1}, \eqref{v21s}, \eqref{v3}, \eqref{v3si}, \eqref{v4} and \eqref{v5} give \eqref{extandagar}. This completes the proof.
\end{proof}
\textbf{Remarks.} 1. Agarwal's identity \eqref{me} can be obtained from \eqref{extandagar} by letting $\e=f=0$ in \eqref{extandagar}, and then applying \eqref{heine} to each of the ${}_3\phi_{2}$'s and to the ${}_2\phi_{1}$ in the third expression on the right side.

2. The fact that the identity \cite[Equation (4.5)]{agarwal} 
\begin{align}\label{agarw}
\sum_{n=0}^{\infty}\frac{(\a)_n}{(\b)_n}t^n=\frac{(\b/\a, q, \a t, q/(\a t); q)_{\infty}}{(q/\a, \b, t, \b/(\a t); q)_{\infty}}+\frac{(1-(q/\b))}{(1-(\a t/\b))}{}_2\phi_{1}(q, q/t; q\b/(\a t);q,  q/\a)
\end{align}
was used in the proof of \eqref{me} (see \cite[p.~294]{agar1}), and \eqref{me} was used in the proof of \eqref{extandagar} given above suggests that a generalization of \eqref{me} for the series 
$\sum_{n=0}^{\infty}\frac{(a_1, a_2, a_3,\cdots, a_r;q)_{n}}{(b_1, b_2, b_3,\cdots, b_r;q)_{n}}t^n$
is not inconceivable.

\noindent
We now prove Theorem \ref{mt1} from Theorem \ref{extvar}. 
\begin{lemma}\label{lemandsim}
If $m$ is a positive integer, we have
\begin{equation}\label{andsim}
\sum_{j=0}^{m}\left[\begin{matrix} m\\j\end{matrix}\right](-a)_j(-a)_{m-j}(-1)^j=\begin{cases}
(q;q^2)_n(a^2;q^2)_n,\hspace{2mm}\text{if}\hspace{2mm} m=2n,\\
0,\hspace{2mm}\text{if}\hspace{2mm} m\hspace{2mm}\text{is odd}.
\end{cases}
\end{equation}
\end{lemma}
\begin{proof}
We note that
\begin{align*}
\sum_{j=0}^{m}\left[\begin{matrix} m\\j\end{matrix}\right](-a)_j(-a)_{m-j}(-1)^j&=(-a)_m\sum_{j=0}^{m}\frac{(q^{-m})_j(-a)_j(q/a)^j}{(q)_{j}(-a^{-1}q^{1-m})_j}\nonumber\\
&=\begin{cases}
(q;q^2)_n(a^2;q^2)_n,\hspace{2mm}\text{if}\hspace{2mm} m=2n,\\
0,\hspace{2mm}\text{if}\hspace{2mm} m\hspace{2mm}\text{is odd},
\end{cases}
\end{align*}
by \cite[p.~526, Equation (1.7)]{gea55}.
\end{proof}
\textbf{Remark.} Ismail and Zhang \cite[Lemma 4.1]{iz} have obtained several interesting results of the similar type as Lemma \ref{lemandsim}.
\begin{proof}[Theorem \textup{\ref{mt1}}][]
Let $\b=-q^2, \g=iq^{3/2}, \d=q^{3/2}, \e=-iq^{3/2}, f=-q^{3/2}, t=q$ in Theorem \ref{extvar}, and then let $\a\to q$. Note that the second expression on the right side of \eqref{extandagar1} vanishes. Hence
{\allowdisplaybreaks\begin{align}\label{mainn1}
&\sum_{n=0}^{\infty}\frac{q^{n}(-q^3;q^2)_n(q;q)_n}{(q^3;q^2)_n(-q^2;q)_n}\nonumber\\
&=\frac{1}{2}\left(1-\frac{1}{q}\right)\frac{(q;q)_{\infty}^{2}}{(-q;q)_{\infty}^{2}}\sum_{n=0}^{\infty}\frac{(-1;q)_n(-q;q^2)_n}{(q;q^2)_n(q;q)_n}q^n\nonumber\\
&\quad+\frac{(-q^3;q^2)_{\infty}}{(q^3;q^2)_{\infty}}\left(1+\frac{1}{q}\right)\sum_{n=0}^{\infty}\frac{(-i\sqrt{q})^n}{(-1)_{n+1}}\sum_{p=0}^{n}\frac{(-1)_p(-q)^p}{(q)_p}\sum_{m=0}^{n}\left[\begin{matrix} n\\m\end{matrix}\right](-i)_m(-i)_{n-m}(-1)^m\nonumber\\
&=\frac{1}{2}\left(1-\frac{1}{q}\right)\frac{(q;q)_{\infty}^{2}}{(-q;q)_{\infty}^{2}}\sum_{n=0}^{\infty}\frac{(-1;q)_n(-q;q^2)_n}{(q;q^2)_n(q;q)_n}q^n\nonumber\\
&\quad+\frac{1}{2q}\frac{(-q;q^2)_{\infty}}{(q^3;q^2)_{\infty}}\sum_{n=0}^{\infty}\frac{(q;q^2)_n(-1;q^2)_n(-q)^n}{(-q)_{2n}}\sum_{p=0}^{2n}\frac{(-1)_p(-q)^p}{(q)_p},
\end{align}}
where in the last step, we applied Lemma \ref{andsim} with $a=i$. This proves \eqref{mainn} upon observing that
\begin{equation}
p_{2n}(-1;q^{-2n-1},-1:q)=\sum_{p=0}^{2n}\frac{(-1)_p(-q)^p}{(q)_p},
\end{equation}
which follows from \eqref{little}. Now \eqref{repre} and \eqref{mainn} imply \eqref{fgf}.
\end{proof}
Next Theorem \ref{congmod4} is proven.

\begin{proof}[Theorem \textup{\ref{congmod4}}][]
Let 
{\allowdisplaybreaks\begin{align}\label{s123}
S_1(q)&:=\sum_{n=0}^{\infty}\frac{q^{n}(-q^3;q^2)_n(q)_n}{(q^3;q^2)_n(-q^2)_n},\nonumber\\
S_2(q)&:=-\frac{1}{2}\left(1-\frac{1}{q}\right)\frac{(q)_{\infty}^{2}}{(-q)_{\infty}^{2}},\nonumber\\
S_3(q)&:=-\frac{1}{2}\left(1-\frac{1}{q}\right)\frac{(q)_{\infty}^{2}}{(-q)_{\infty}^{2}}\sum_{n=1}^{\infty}\frac{(-1)_n(-q;q^2)_n}{(q;q^2)_n(q)_n}q^n.
\end{align}}
By \eqref{mainn}, proving \eqref{mod4} is equivalent to showing
\begin{equation}\label{mod4a}
S_1(q)+S_2(q)+S_3(q)\equiv\frac{1}{2q}-\frac{1}{2}\hspace{1mm}(\textup{mod}\hspace{1mm}4).
\end{equation}
Note that $(q)_{\infty}^{2} \equiv {(-q)_{\infty}^{2}} \hspace{1mm}(\textup{mod}\hspace{1mm}4)$. Hence
\begin{align}\label{s13}
S_1(q)+S_3(q)&\equiv S_1(q)+\sum_{n=1}^{\infty}\frac{(-q)_{n-1}(-q;q^2)_nq^{n-1}}{(q^3;q^2)_{n-1}(q)_n}
\hspace{1mm}(\textup{mod}\hspace{1mm}4)\nonumber\\
&=S_1(q)+\sum_{n=0}^{\infty}\frac{(-q)_{n}(-q;q^2)_{n+1}q^{n}}{(q^3;q^2)_{n}(q)_{n+1}}\nonumber\\
&=\sum_{n=0}^{\infty}\frac{q^{n}(-q;q^2)_{n+1}(q)_n}{(q^3;q^2)_n(-q)_{n+1}}+\sum_{n=0}^{\infty}\frac{(-q)_{n}(-q;q^2)_{n+1}q^{n}}{(q^3;q^2)_{n}(q)_{n+1}}\nonumber\\
&=\sum_{n=0}^{\infty}\frac{(-q;q^2)_{n+1}q^n}{(q^3;q^2)_n}\left(\frac{(q)_n}{(-q)_{n+1}}+\frac{(-q)_{n}}{(q)_{n+1}}\right)\nonumber\\
&=\sum_{n=0}^{\infty}\frac{(-q;q^2)_{n+1}q^n}{(q^2;q)_{2n+1}}\left((q)_{n}^2(1-q^{n+1})+(-q)_{n}^{2}(1+q^{n+1})\right)\nonumber\\
&\equiv 2\sum_{n=0}^{\infty}\frac{(-q;q^2)_{n+1}(q)_{n}^2q^n}{(q^2;q)_{2n+1}}\hspace{1mm}(\textup{mod}\hspace{1mm}4),
\end{align}
since $(q)_{n}^{2}\equiv (-q)_{n}^{2}\hspace{1mm}(\textup{mod}\hspace{1mm}4)$. Now
{\allowdisplaybreaks\begin{align}\label{s2}
S_2(q)&=-\frac{1}{2}\left(1-\frac{1}{q}\right)\left(1+2\sum_{n=1}^{\infty}(-1)^nq^{n^2}\right)^2\nonumber\\
&=\frac{1}{2q}-\frac{1}{2}+\frac{2}{q}(1-q)\left(\sum_{n=1}^{\infty}(-1)^nq^{n^2}+\left(\sum_{n=1}^{\infty}(-1)^nq^{n^2}\right)^2\right).
\end{align}}
From \eqref{s13} and \eqref{s2}, it suffices to show that
\begin{equation}\label{mod4b}
2\sum_{n=0}^{\infty}\frac{(-q;q^2)_{n+1}(q)_{n}^2q^n}{(q^2;q)_{2n+1}}\equiv-\frac{2}{q}(1-q)\left(\sum_{n=1}^{\infty}(-1)^nq^{n^2}+\left(\sum_{n=1}^{\infty}(-1)^nq^{n^2}\right)^2\right)\hspace{1mm}(\textup{mod}\hspace{1mm}4),
\end{equation}
or equivalently, 
\begin{equation}\label{mod4c}
\sum_{n=0}^{\infty}\frac{(q^3;q^2)_{n}(q^2;q^2)_{n}q^{n+1}}{(q^2;q^2)_{n+1}(q^3;q^2)_n}\equiv-\left(\sum_{n=1}^{\infty}(-1)^nq^{n^2}+\left(\sum_{n=1}^{\infty}(-1)^nq^{n^2}\right)^2\right)\hspace{1mm}(\textup{mod}\hspace{1mm}2).
\end{equation}
Now
\begin{align}
\sum_{n=0}^{\infty}\frac{(q^3;q^2)_{n}(q^2;q^2)_{n}q^{n+1}}{(q^2;q^2)_{n+1}(q^3;q^2)_n}&=\sum_{n=0}^{\infty}\frac{q^{n+1}}{1-q^{2n+2}}\nonumber\\
&=\sum_{N=1}^{\infty}d_{o}(N)q^{N},
\end{align}
where $d_{o}(N)$ is the number of odd divisors of $N$. Also,
\begin{align}
-\left(\sum_{n=1}^{\infty}(-1)^nq^{n^2}+\left(\sum_{n=1}^{\infty}(-1)^nq^{n^2}\right)^2\right)\equiv\sum_{n=1}^{\infty}q^{n^2}+\left(\sum_{n=1}^{\infty}q^{n^2}\right)^2\hspace{1mm}(\textup{mod}\hspace{1mm}2).
\end{align}
Let
\begin{equation}
\sum_{N=1}^{\infty}a(N)q^{N}:=\sum_{n=1}^{\infty}q^{n^2}+\left(\sum_{n=1}^{\infty}q^{n^2}\right)^2.
\end{equation}
Let $r_2(m)$ denote the number of representations of $m$ as a sum of two squares, where representations with different orders or different signs of the
summands are regarded as distinct. Now if $N$ is not a square, then the number of representations of $N$ as a sum of $2$ positive squares is equal to 
\begin{align}
\tfrac{1}{4}r_{2}(N)=d_1(N)-d_3(N)\equiv d_{o}(N)\hspace{1mm}(\textup{mod}\hspace{1mm}2),
\end{align}
where Jacobi's formula was employed in the penultimate step. That is, $a(N)\equiv d_{o}(N)\hspace{1mm}(\textup{mod}\hspace{1mm}2)$.
If, however, $N$ is a square, then the number of representations of $N$ as a sum of $2$ positive squares is equal to $\frac{1}{4}r_2(n)-1$. Hence,
\begin{align}
a(N)=\frac{1}{4}r_{2}(N)=d_1(N)-d_3(N)\equiv d_{o}(N)\hspace{1mm}(\textup{mod}\hspace{1mm}2).
\end{align}
This implies that \eqref{mod4c} always holds, and this proves the theorem.
\end{proof}

\section{Congruences for $\overline{p}_{\omega}(n)$}\label{pwcong}
This section is devoted to proving Theorem \ref{pw4}. We start with the following series $S(q)$:
\begin{align}
S(q):=\sum_{n=1}^{\infty} \frac{q^n (q^{n+1};q)_n (q^{2n+2};q^2)_{\infty}}{(1+q^n)(-q^{n+1};q)_n (-q^{2n+2};q^2)_{\infty}}. \label{S}
\end{align}
\begin{lemma}\label{Sq}
The following identity holds:
\begin{align}
S(q)=-\sum_{n=1}^{\infty} (-1)^n q^{2n^2}
+  \sum_{k=1}^{\infty} \frac{q^{k} (q;q^2)_{k}}{(-q;q^2)_k (1+q^{2k})}. \label{s(q)}
\end{align}
\end{lemma}
\proof
\allowdisplaybreaks{
\begin{align*}
S(q)=&
\sum_{n=1}^{\infty} \frac{q^n (q^{n+1};q)_n (q^{2n+2};q^2)_{\infty}}{(1+q^n)(-q^{n+1};q)_n (-q^{2n+2};q^2)_{\infty}}&
\\& =\frac{(q;q)_{\infty}}{(-q;q)_{\infty}} \sum_{n=1}^{\infty} \frac{q^n (-q;q)_{n-1}(-q^{2n+1};q^2)_{\infty}}{(q;q)_n (q^{2n+1};q^2)_{\infty}}\\
&=\frac{(q;q)_{\infty}}{(-q;q)_{\infty}} \sum_{n=1}^{\infty} \frac{ q^n (-q;q)_{n-1}}{(q;q)_n } \sum_{k=0}^{\infty} \frac{(-1;q^2)_{k}}{(q^2;q^2)_k} q^{(2n+1)k}\\
&=\frac{(q;q)_{\infty}}{(-q;q)_{\infty}} \sum_{k=0}^{\infty} \frac{q^{k} (-1;q^2)_{k}}{(q^2;q^2)_k} \sum_{n=1}^{\infty} \frac{ (-q;q)_{n-1}}{(q;q)_n } q^{(2k+1)n} \\
&=\frac{(q;q)_{\infty}}{(-q;q)_{\infty}} \sum_{k=0}^{\infty} \frac{q^{k} (-1;q^2)_{k}}{(q^2;q^2)_k} \left(-\frac{1}{2}+ \frac{1}{2}\sum_{n=0}^{\infty} \frac{ (-1;q)_{n}}{(q;q)_n } q^{(2k+1)n}\right) \\
&=-\frac{(q;q)_{\infty}}{2(-q;q)_{\infty}} \sum_{k=0}^{\infty} \frac{q^{k} (-1;q^2)_{k}}{(q^2;q^2)_k} + \frac{(q;q)_{\infty}}{2(-q;q)_{\infty}} \sum_{k=0}^{\infty} \frac{q^{k} (-1;q^2)_{k}}{(q^2;q^2)_k} \sum_{n=0}^{\infty} \frac{ (-1;q)_{n}}{(q;q)_n } q^{(2k+1)n} \\
&=-\frac{(q;q)_{\infty}}{2(-q;q)_{\infty}}\frac{(-q;q^2)_\infty}{(q;q^2)_{\infty}}
+ \frac{(q;q)_{\infty}}{2(-q;q)_{\infty}}\sum_{k=0}^{\infty} \frac{q^{k} (-1;q^2)_{k}}{(q^2;q^2)_k} \frac{(-q^{2k+1};q)_{\infty}}{(q^{2k+1};q)_\infty } \\
&=-\frac{(q^2;q^2)_{\infty}}{2(-q^2;q^2)_{\infty}}
+ \frac{1}{2}\sum_{k=0}^{\infty} \frac{q^{k} (-1;q^2)_{k}}{(q^2;q^2)_k} \frac{(q;q)_{2k}}{(-q;q)_{2k}} \\
&=-\frac{1}{2}-\sum_{n=1}^{\infty} (-1)^n q^{2n^2}
+ \frac{1}{2}+ \sum_{k=1}^{\infty} \frac{q^{k} (q;q^2)_{k}}{(-q;q^2)_k (1+q^{2k})},
\end{align*}
}
where \eqref{qbin} was used for the third and seventh equalities. 
\endproof

Let
\begin{equation}
A(q):=
\sum_{k=1}^{\infty} \frac{q^{k} (q;q^2)_{k}}{(-q;q^2)_k (1+q^{2k})}. \label{Aq}
\end{equation}

\begin{lemma} \label{prop2}
we have
\begin{equation*}
A(q)+A(-q)=
-\frac{1}{2}+\frac{1}{2} \frac{(q^2,q^2;q^2)_{\infty}}{(-q^2,-q^2;q^2)_{\infty}}.
\end{equation*}
\end{lemma}
\proof
We now extend the sum in \eqref{Aq} to negative infinity:
\begin{align}
\sum_{k=-\infty}^{\infty} \frac{q^{k} (q;q^2)_{k}}{(-q;q^2)_k (1+q^{2k})} =\frac{1}{2}+\sum_{k=1}^{\infty}  \frac{q^{k} (q;q^2)_{k}}{(-q;q^2)_k (1+q^{2k})} +\sum_{k=1}^{\infty}  \frac{(-1)^k q^{k} (-q;q^2)_{k}}{(q;q^2)_k (1+q^{2k})}, \label{7}
\end{align}
where $(a;q)_{n}:=(a;q)_{\infty}/(aq^n;q)_{\infty}$. 
Thus
\begin{align}
\sum_{k=-\infty}^{\infty} \frac{q^{k} (q;q^2)_{k}}{(-q;q^2)_k (1+q^{2k})}=\frac{1}{2}+A(q)+A(-q). 
\label{8}
\end{align}
Set $q\to q^2$, $a=-1, c=q, d=1, e=q$, and $f=-1$ in \eqref{9}. 
Then we obtain
\allowdisplaybreaks{
\begin{align*}
\sum_{k=-\infty}^{\infty} \frac{q^{k} (q;q^2)_{k}}{(-q;q^2)_k (1+q^{2k})}&=
\frac{1}{2}\sum_{k=-\infty}^{\infty} \frac{q^{k} (q;q^2)_{k} (-1;q^2)_k}{(-q;q^2)_k (-q^2;q^2)_k}\\
&=\frac{1}{2}\frac{(q,q^2,-q,q^2;q^2)_{\infty}}{(-q^2, -q^2,-q, q;q^2)_{\infty}}\\
&=\frac{1}{2} \frac{(q^2,q^2;q^2)_{\infty}}{(-q^2,-q^2;q^2)_{\infty}},
\end{align*}}
which with \eqref{8} completes the proof.
\endproof

\begin{lemma} \label{lem3}
We have
\begin{align}
& \frac{1}{2}\bigg( A(q)-A(-q) \bigg) \equiv  q\frac{(q^8;q^8)_{\infty}^4}{(q^4;q^4)_{\infty}^2} \pmod{4}. \label{63}
\end{align}
\end{lemma}
\proof
First note that using Alladi's identity \cite[p.~215]{alladi}
we obtain
\begin{align*}
\frac{(q;q^2)_k}{(-q;q^2)_k}=1-2\sum_{j=1}^k \frac{q^{2j-1} (q;q^2)_{j-1}}{(-q;q^2)_{j}}.
\end{align*}
Thus
\begin{align}
A(q)&=\sum_{k=1}^{\infty} \frac{q^{k} (q;q^2)_{k}}{(-q;q^2)_k (1+q^{2k})} \notag\\
&=\sum_{k=1}^{\infty} \frac{q^{k}}{1+q^{2k}}-2\sum_{k=1}^{\infty} \frac{q^{k}}{1+q^{2k}} \sum_{j=1}^k \frac{q^{2j-1} (q;q^2)_{j-1}}{(-q;q^2)_{j}} \notag\\
&\equiv  \sum_{k=1}^{\infty} \frac{q^{k}}{1+q^{2k}}+ 2\sum_{k=1}^{\infty} \frac{q^{k}}{1+q^{2k}} \sum_{j=1}^k \frac{q^{2j-1}}{1+q^{2j-1}} \pmod{4}. \label{333}
\end{align}
Now
\begin{align}
\sum_{k=1}^{\infty} \frac{q^{k}}{1+q^{2k}} \sum_{j=1}^k \frac{q^{2j-1}}{1+q^{2j-1}}& = \sum_{k=1}^{\infty} \frac{q^{k}}{1+q^{2k}} \sum_{j=1}^k \sum_{m=1}^{\infty} (-1)^{m-1} q^{(2j-1)m}  \notag\\
&= \sum_{k=1}^{\infty} \frac{q^{k}}{1+q^{2k}} \sum_{m=1}^{\infty} \sum_{j=1}^k (-1)^{m-1} q^{(2j-1)m} \notag\\
&= -\sum_{k=1}^{\infty} \frac{q^{k}}{1+q^{2k}} \sum_{m=1}^{\infty} (-q)^m \sum_{j=0}^{k-1}  q^{2mj} \notag\\
&= - \sum_{k=1}^{\infty} \frac{q^{k}}{1+q^{2k}} \sum_{m=1}^{\infty} (-q)^m \frac{(1-q^{2km})}{1-q^{2m}} \notag \\
&\equiv \sum_{k,m=1}^{\infty} \frac{q^{k+m}}{(1+q^{2k})(1+q^{2m})} + \sum_{k,m=1}^{\infty}  \frac{q^{2km+k+m}}{(1+q^{2k})(1+q^{2m})} \pmod{2} \notag\\
&\equiv \sum_{k=1}^{\infty} \frac{q^{2k}}{(1-q^{2k})^2}+\sum_{k=1}^{\infty} \frac{q^{2k^2+2k}}{(1-q^{2k})^2} \pmod{2}, \label{334}
\end{align}%
where the last congruence follows since each of the double summations is symmetric in $k$ and $m$. Thus, by \eqref{333} and \eqref{334}
\begin{align*}
A(q)
&\equiv  \sum_{k=1}^{\infty} \frac{q^{k}}{1+q^{2k}}+ 2\sum_{k=1}^{\infty} \frac{q^{2k}}{1-q^{4k}} +2 \sum_{k=1}^{\infty} \frac{q^{2k(k+1)}}{1-q^{4k}} \pmod{4}\\
&\equiv \sum_{k=1}^{\infty} \frac{q^{2k-1}}{1+q^{4k-2}}+\sum_{k=1}^{\infty} \frac{q^{2k}}{1+q^{4k}}+ 2\sum_{k=1}^{\infty} \frac{q^{2k}}{1-q^{4k}} +2\sum_{k=1}^{\infty} \frac{q^{2k(k+1)}}{1-q^{4k}} \pmod{4} .
\end{align*}
Since the odd powers of $q$ appear only in the first sum on the right hand side above, we see that
\begin{align*}
\frac{1}{2}\bigg(A(q)-A(-q)\bigg) & \equiv \sum_{k=1}^{\infty} \frac{q^{2k-1}}{1+q^{4k-2}} \pmod{4}\\
&=q\frac{(q^8;q^8)_{\infty}^4}{(q^4;q^4)_{\infty}^2},
\end{align*}
where the last equality follows from (32.26) in \cite{fine}. It can also be derived by letting $q\to q^4$, and then substituting $a=-q^{-2}, b=-q^2, z=q^2$ in Ramanujan's $_1 \psi_1$ summation formula \cite[p.239, (II 29)]{gasper} 
\begin{equation}\label{1psi1sf}
\sum_{n=-\infty}^{\infty}\frac{(a;q)_{n}}{(b;q)_{n}}z^n=\frac{(az;q)_{\infty}(q/(az);q)_{\infty}(q;q)_{\infty}(b/a;q)_{\infty}}{(z;q)_{\infty}(b/(az);q)_{\infty}(b;q)_{\infty}(q/a;q)_{\infty}},
\end{equation}
valid for $|b/a|<|z|<1$ and $|q|<1$. 
\endproof

By Lemmas~\ref{prop2} and \ref{lem3}, we have
\begin{equation*}
A(q)\equiv -\frac{1}{4}+ \frac{1}{4}\frac{(q^2;q^2)_{\infty}^2}{(-q^2;q^2)_{\infty}^2}+ q\frac{(q^8;q^8)_{\infty}^4}{(q^4;q^4)_{\infty}^2} \pmod{4}. 
\end{equation*}
Also, we recall that
\begin{align}
\phi(-q)&:=\sum_{n=-\infty}^{\infty}  (-1)^n q^{n^2}=(q;q)_{\infty}(q;q^2)_{\infty}=\frac{(q;q)_{\infty}}{(-q;q)_{\infty}}, \label{phi}\\
\psi(q)&:=\sum_{n=0}^{\infty} q^{n(n+1)/2}=\frac{(q^2;q^2)_{\infty}}{(q;q^2)_{\infty}}=\frac{(q^2;q^2)^2_{\infty}}{(q;q)_{\infty}}. \label{psi}
\end{align}
Thus,
\begin{align}
A(q)\equiv -\frac{1}{4}+\frac{1}{4}\sum_{m,n=-\infty}^{\infty} (-1)^{m+n} q^{2(m^2+n^2)} + q \sum_{m,n=0}^{\infty} q^{2m(m+1)+2n(n+1)} \pmod{4}. \label{a(q)}
\end{align}
We are now ready to prove Theorem~\ref{pw4}.
  First, note that
\begin{align}\label{4mod}
\frac{1+x}{1-x} \equiv \frac{1-x}{1+x} \pmod{4}. 
\end{align}
Thus from \eqref{gfbpo},
\begin{align}
\sum_{n=1}^{\infty} \overline{p}_{\omega}(n)q^n
&\equiv\sum_{n=1}^{\infty} \frac{q^n (q^{n+1};q)_n (q^{2n+2};q^2)_{\infty}}{(1-q^n)(-q^{n+1};q)_n (-q^{2n+2};q^2)_{\infty}} \pmod{4}. \label{22}
\end{align}
Now
\begin{align}
&\sum_{n=1}^{\infty} \frac{q^n (q^{n+1};q)_n (q^{2n+2};q^2)_{\infty}}{(1-q^n)(-q^{n+1};q)_n (-q^{2n+2};q^2)_{\infty}} \notag\\
&=\sum_{n=1}^{\infty} \frac{q^n (q^{n+1};q)_n (q^{2n+2};q^2)_{\infty}}{(1+q^n)(-q^{n+1};q)_n (-q^{2n+2};q^2)_{\infty}} + 2 \sum_{n=1}^{\infty} \frac{q^{2n} (q^{n+1};q)_n (q^{2n+2};q^2)_{\infty}}{(1-q^{2n})(-q^{n+1};q)_n (-q^{2n+2};q^2)_{\infty}} \notag \\
&\equiv \sum_{n=1}^{\infty} \frac{q^n (q^{n+1};q)_n (q^{2n+2};q^2)_{\infty}}{(1+q^n)(-q^{n+1};q)_n (-q^{2n+2};q^2)_{\infty}}+2 \sum_{n=1}^{\infty} \frac{q^{2n} }{1-q^{2n}} \pmod{4} \notag\\
&\equiv \sum_{n=1}^{\infty} \frac{q^n (q^{n+1};q)_n (q^{2n+2};q^2)_{\infty}}{(1+q^n)(-q^{n+1};q)_n (-q^{2n+2};q^2)_{\infty}} +2 \sum_{n=1}^{\infty} q^{2n^2}, \label{3}
\end{align}
where the second last congruence follows from the fact that 
\begin{equation}\label{2mod}
1+x\equiv 1-x \pmod{2}.
\end{equation}
For the last congruence above, we used Clausen's identity \cite[p.~16, Equation (14.51)]{fine} 
\begin{equation*}
\sum_{n=1}^{\infty}d(n)q^n=\sum_{n=1}^{\infty}\left(\frac{1+q^n}{1-q^n}\right)q^{n^2},
\end{equation*}
which implies that
\begin{align}
\sum_{n=1}^{\infty} \frac{q^{n} }{1-q^{n}}
\equiv \sum_{n=1}^{\infty} q^{ n^2} \pmod{2}. \label{44}
\end{align}
Thus, from \eqref{S}, \eqref{22} and \eqref{3}, we have
\begin{align}
\sum_{n=1}^{\infty} \overline{p}_{\omega}(n)q^n&\equiv  S(q) +2\sum_{n=1}^{\infty} q^{2n^2} \pmod{4}. \label{thm4.4}
\end{align}

\begin{theorem} \label{conj}
We have
\begin{align}
\sum_{n=1}^{\infty} \overline{p}_{\omega}(n)q^n \equiv -1 +\sum_{m,n=0}^{\infty} (-1)^{m+n} q^{2(m^2+n^2)} + q \sum_{m,n=0}^{\infty} q^{2m(m+1)+2n(n+1)} \pmod{4}. 
\end{align}
\end{theorem}
\proof
The congruence follows from  \eqref{s(q)}, \eqref{Aq}, \eqref{a(q)}, and \eqref{thm4.4}. 
\endproof

\noindent
Thus it immediately follows from Theorem~\ref{conj} that $\overline{p}_{\omega}(4n+3)\equiv 0 \pmod{4}$. Since $8n+6=2(4n+3)$ and $4n+3$ cannot be written as a sum of two squares, this also proves $\overline{p}_{\omega}(8n+6)\equiv 0 \pmod{4}$.
\section{Different representations of the generating function of $\overline{\textup{spt}}_{\omega}(n)$}\label{sptwdf}

As in \cite{jennings-shaffer}, define
\begin{align*}
\sum_{n=1}^{\infty} \overline{\textup{spt}}(n) q^n =\sum_{n=1}^{\infty} \frac{q^n (-q^{n+1};q)_{\infty}}{(1-q^n)^2(q^{n+1};q)_{\infty}}.
\end{align*}
Here we note that this $\overline{\textup{spt}}(n)$ is exactly a half of the smallest parts function for overpartitions of $n$ defined by Bringmann, Lovejoy, and Osburn in \cite{bringmannlovejoyosburn1}. By taking $d=1$ and $e=0$ in 
Equations (1.1), (1.2), and Theorem 7.1. in \cite{bringmannlovejoyosburn2}, we see that
\begin{align}
\sum_{n=1}^{\infty} \overline{\textup{spt}}(n) q^n =\sum_{n=1}^{\infty} \frac{q^n (-q^{n+1};q)_{\infty}}{(1-q^n)^2(q^{n+1};q)_{\infty}} =\frac{(-q;q)_{\infty}}{(q;q)_{\infty}} \sum_{n=1}^{\infty} \frac{nq^n}{1-q^n} + 2\frac{(-q;q)_{\infty}}{(q;q)_{\infty}}\sum_{n=1}^{\infty} \frac{(-1)^n q^{n^2+n}}{(1-q^n)^2}. \label{barspt}
\end{align} 
We also define
\begin{align*}
\sum_{n=1}^{\infty} \overline{\textup{spt}}_2(n)q^n=\sum_{n=1}^{\infty} \frac{q^{2n}(-q^{2n+1};q)_{\infty}}{(1-q^{2n})^2(q^{2n+1};q)_{\infty}}.
\end{align*}
Then, again by taking $d=1, e=q^{-1}$ and replacing $q$ by $q^2$ in Equations (1.1), (1.2), and Theorem 7.1. in \cite{bringmannlovejoyosburn2}, we have
\begin{align}
\sum_{n=1}^{\infty} \overline{\textup{spt}}_2(n)q^n = \sum_{n=1}^{\infty} \frac{q^{2n}(-q^{2n+1};q)_{\infty}}{(1-q^{2n})^2(q^{2n+1};q)_{\infty}}=\frac{(-q;q)_{\infty}}{(q;q)_{\infty}}\sum_{n=1}^{\infty} \frac{nq^{n}}{1-q^{n}}+\frac{(-q;q)_{\infty}}{(q;q)_{\infty}}\sum_{n=1}^{\infty} \frac{(-1)^n q^n (1+q^{n^2})}{(1-q^n)^2}. \label{barspt2}
\end{align}

\begin{theorem} \label{thmbarsptomega}
We have
\begin{align}
\sum_{n=1}^{\infty} \overline{\textup{spt}}_{\omega}(n)q^n &=\frac{(-q^2;q^2)_{\infty}}{(q^2;q^2)_{\infty}} \sum_{n=1}^{\infty} \frac{nq^{n}}{1-q^n} +2 \frac{(-q^2;q^2)_{\infty} }{(q^2;q^2)_{\infty}} \sum_{n=1}^{\infty} \frac{(-1)^n q^{2n(n+1)}}{(1-q^{2n})^2}. \label{barsptomega}
\end{align}
\end{theorem}
\proof
In \eqref{baitra}, we set $a=1$, $p_1=z=p_2^{-1}$, and $f=-1$. Then we obtain
\begin{align}
\sum_{n=0}^{\infty} \frac{(z;q)_n(z^{-1};q)_n (-q;q^2)_n q^n}{(q;q)_n (q;q^2)_n (-q;q)_n} =\frac{(zq;q)_{\infty} (z^{-1}q;q)_{\infty}}{(q;q)_{\infty}^2}\left(1+2 \sum_{n=1}^{\infty} \frac{(1-z)(1-z^{-1}) (-1)^n  q^{2n(n+1)}}{(1-zq^{2n})(1-z^{-1}q^{2n})}\right). \label{bailey}
\end{align}
Now take the second derivative on both sides of \eqref{bailey} with respect to $z$, then set $z=1$ to obtain
\begin{align*}
\sum_{n=1}^{\infty} \frac{q^n (q;q)_n (-q;q^2)_n }{(1-q^n)^2 (-q;q)_n(q;q^2)_{n} } =\sum_{n=1}^{\infty} \frac{nq^n}{1-q^n} +2 \sum_{n=1}^{\infty} \frac{ (-1)^n q^{2n(n+1)}}{(1-q^{2n})^2}.
\end{align*}
Multiply both sides of the above identity by $(-q^2;q^2)_{\infty}/(q^2;q^2)_{\infty}$ to get
\begin{align}
\frac{(-q^2;q^2)_{\infty}}{(q^2;q^2)_{\infty}}\sum_{n=1}^{\infty} \frac{q^n (q;q)_n (-q;q^2)_n }{(1-q^n)^2 (-q;q)_n(q;q^2)_{n} }  =\frac{(-q^2;q^2)_{\infty}}{(q^2;q^2)_{\infty}} \sum_{n=1}^{\infty} \frac{nq^{n}}{1-q^n} +2 \frac{(-q^2;q^2)_{\infty} }{(q^2;q^2)_{\infty}} \sum_{n=1}^{\infty} \frac{(-1)^n q^{2n(n+1)}}{(1-q^{2n})^2}. \label{2}
\end{align}
Note that the left hand side of \eqref{2} is 
\begin{align}
\frac{(-q^2;q^2)_{\infty}}{(q^2;q^2)_{\infty}} \sum_{n=1}^{\infty} \frac{q^n (q;q)_n (-q;q^2)_n }{(1-q^n)^2 (-q;q)_n(q;q^2)_{n} }&=\frac{(-q^2;q^2)_{\infty}}{(q^2;q^2)_{\infty}} \sum_{n=1}^{\infty} \frac{q^n (q;q)_n (-q;q)_{2n} (q^2;q^2)_n }{(1-q^n)^2 (-q;q)_n(q;q)_{2n} (-q^2;q^2)_n} \notag\\
&=\frac{(-q^2;q^2)_{\infty}}{(q^2;q^2)_{\infty}} \sum_{n=1}^{\infty} \frac{q^n (-q^{n+1};q)_{n} (q^2;q^2)_n }{(1-q^n)^2 (q^{n+1};q)_{n} (-q^2;q^2)_n} \notag\\
&=\sum_{n=1}^{\infty} \frac{q^n (-q^{n+1};q)_{n} (-q^{2n+2};q^2)_\infty }{(1-q^n)^2 (q^{n+1};q)_{n} (q^{2n+2};q^2)_\infty}\notag \\
&=\sum_{n=1}^{\infty} \overline{\textup{spt}}_{\omega}(n)q^n, 
\notag 
\end{align}
where the last equality follows from the definition of $\overline{\textup{spt}}_{\omega}(n)$ in \eqref{sptbaromega}.  This completes the proof. 
\endproof


\begin{lemma} \label{lem1}
We have
\begin{align*}
\sum_{n=1}^{\infty} \frac{(-1)^n q^{n}}{(1-q^n)^2}=-\sum_{n=1}^{\infty} \frac{(2n-1)q^{2n-1}}{1-q^{2n-1}}.
\end{align*}
\end{lemma}

\proof
Note that
\begin{align*}
\sum_{n=1}^{\infty} \frac{(-1)^n q^{n}}{(1-q^n)^2}=\sum_{n,k=1}^{\infty}  (-1)^n k q^{kn} =-\sum_{k=1}^{\infty} \frac{kq^{k}}{1+q^{k}}.
\end{align*}
Robbins \cite[Theorem 3]{robbins} has shown that this is the negative of the generating function for the sum of the odd divisors of $n$ which is the right-hand side of the identity we wanted to prove.
\endproof

We now relate our $\overline{\textup{spt}}_{\omega}(n)$ to $\overline{\textup{spt}}(n)$ and $\overline{\textup{spt}}_2(n)$.  
\begin{corollary}
We have
\begin{align}
\sum_{n=1}^{\infty} \overline{\textup{spt}}_{\omega}(n)q^n 
&=\frac{(-q^2;q^2)_{\infty}}{(q^2;q^2)_{\infty}} \left( \sum_{n=1}^{\infty} \frac{(2n-1)q^{2n-1}}{1-q^{2n-1}}+ \sum_{n=1}^{\infty} \frac{n q^{2n}}{1-q^{2n}} \right)+ \sum_{n=1}^{\infty} \overline{\textup{spt}}(n)q^{2n} , \label{4}\\
&= \frac{(-q^2;q^2)_{\infty}}{(q^2;q^2)_{\infty}} \left( \sum_{n=1}^{\infty} \frac{(2n-1)q^{2n-1}}{1-q^{2n-1}} +2\sum_{n=1}^{\infty} \frac{(2n-1) q^{4n-2}}{1-q^{4n-2}} \right)+ 2\sum_{n=1}^{\infty} \overline{\textup{spt}}_2(n)q^{2n}. \label{40}
\end{align}
\end{corollary}

\proof From \eqref{barsptomega},
{\allowdisplaybreaks\begin{align*}
\sum_{n=1}^{\infty} \overline{\textup{spt}}_{\omega}(n)q^n 
&=\frac{(-q^2;q^2)_{\infty}}{(q^2;q^2)_{\infty}} \sum_{n=1}^{\infty} \frac{(2n-1)q^{2n-1}}{1-q^{2n-1}}+\frac{(-q^2;q^2)_{\infty}}{(q^2;q^2)_{\infty}} \sum_{n=1}^{\infty} \frac{2n q^{2n}}{1-q^{2n}} + 2 \frac{(-q^2;q^2)_{\infty} }{(q^2;q^2)_{\infty}} \sum_{n=1}^{\infty} \frac{(-1)^n q^{2n(n+1)}}{(1-q^{2n})^2}\notag\\
&=\frac{(-q^2;q^2)_{\infty}}{(q^2;q^2)_{\infty}} \sum_{n=1}^{\infty} \frac{(2n-1)q^{2n-1}}{1-q^{2n-1}}+ \frac{(-q^2;q^2)_{\infty} }{(q^2;q^2)_{\infty}} \sum_{n=1}^{\infty} \frac{n q^{2n}}{1-q^{2n}}+  \sum_{n=1}^{\infty} \overline{\textup{spt}}(n)q^{2n} \notag\\
&= \frac{(-q^2;q^2)_{\infty}}{(q^2;q^2)_{\infty}} \left( \sum_{n=1}^{\infty} \frac{(2n-1)q^{2n-1}}{1-q^{2n-1}}+ \sum_{n=1}^{\infty} \frac{n q^{2n}}{1-q^{2n}} \right)+ \sum_{n=1}^{\infty} \overline{\textup{spt}}(n)q^{2n} , 
\end{align*}}
where the second last equality follows from \eqref{barspt}. Also, by \eqref{barspt2},
\begin{align*}
\sum_{n=1}^{\infty} \overline{\textup{spt}}_{\omega}(n)q^n 
&= \frac{(-q^2;q^2)_{\infty}}{(q^2;q^2)_{\infty}} \left( \sum_{n=1}^{\infty} \frac{(2n-1)q^{2n-1}}{1-q^{2n-1}}- 2\sum_{n=1}^{\infty} \frac{(-1)^n q^{2n}}{(1-q^{2n})^2} \right)+ 2\sum_{n=1}^{\infty} \overline{\textup{spt}}_2(n)q^{2n}, 
\end{align*}
which with Lemma~\ref{lem1} yields \eqref{40}. 
\endproof

\section{Congruences for $\overline{\textup{spt}}_{\omega}(n)$}\label{sptwcong}
The congruences satisfied by $\overline{\textup{spt}}_{\omega}(n)$, which are given in Theorem \ref{thm1}, are proved in this section.

\subsection{Congruences modulo $3$}\label{cong3}
We prove \eqref{equ_3n0} and \eqref{equ_3n2} here. Let 
\begin{align}
S(q):=\sum_{n=1}^{\infty} c_nq^n=\frac{(-q^2;q^2)_{\infty}}{(q^2;q^2)_{\infty}}\left(\sum_{n=1}^{\infty} \frac{(2n-1)q^{2n-1}}{1-q^{2n-1}} +2\sum_{n=1}^{\infty} \frac{(2n-1) q^{4n-2}}{1-q^{4n-2}}\right). \label{s22}
\end{align}
Then, by \eqref{40},
\begin{align}\label{extra}
\overline{\textup{spt}}_{\omega}(n)=c_n + 2 \overline{\textup{spt}}_2\left(\frac{n}{2}\right),
\end{align}
where we follow the convention that $\overline{\textup{spt}}_2(x)=0$ if $x$ is not a positive integer. 

By \eqref{eq3n0} and \eqref{eq3n1}, it suffices to show that $c_{3n}\equiv c_{3n+2}\equiv 0 \pmod{3}$. Now
\begin{align*}
S(q) 
&\equiv \frac{(-q^2;q^2)_{\infty}}{(q^2;q^2)_{\infty}} \left( \sum_{n=1}^{\infty} \frac{(2n-1)q^{2n-1}}{1-q^{2n-1}} - \sum_{n=1}^{\infty} \frac{(2n-1)q^{4n-2}}{1-q^{4n-2}}  \right) \pmod{3}\\
&=\frac{(-q^2;q^2)_{\infty}}{(q^2;q^2)_{\infty}} \sum_{n=1}^{\infty} \frac{(2n-1)q^{2n-1}}{1-q^{4n-2}} \\
&=\frac{(-q^2;q^2)_{\infty}}{(q^2;q^2)_{\infty}}\frac{q(q^4;q^4)_{\infty}^8}{(q^2;q^2)_{\infty}^4} \\
&=q(-q^2;q^2)_{\infty}^9(q^2;q^2)_{\infty}^3\\
&\equiv q (-q^6;q^6)_{\infty}^3 (q^6;q^6)_{\infty} \pmod{3}, 
\end{align*}
where the third equality follows from  \cite[Equation (32.31)]{fine}. Hence $c_{3n}\equiv c_{3n+2}\equiv 0 \pmod{3}$.


\subsection{Another proof of \eqref{equ_3n0}} 
Let 
\begin{align*}
M_1(q):=\sum_{n=1}^{\infty} d_nq^n=\frac{(-q^2;q^2)_{\infty}}{(q^2;q^2)_{\infty}} \left( \sum_{n=1}^{\infty} \frac{(2n-1)q^{2n-1}}{1-q^{2n-1}}+ \sum_{n=1}^{\infty} \frac{n q^{2n}}{1-q^{2n}} \right).
\end{align*}
Since \cite[Thm. 1.2]{jennings-shaffer} implies $\overline{\textup{spt}}(3n)\equiv 0 \pmod{3}$, it suffices to show $d_{3n}\equiv 0 \pmod{3}$ by \eqref{4}. 
Now
\begin{align*}
M_1(q)&=\frac{(-q^2;q^2)_{\infty}}{(q^2;q^2)_{\infty}} \left( \sum_{n=1}^{\infty} \frac{(2n-1)q^{2n-1}}{1-q^{2n-1}}+\sum_{n=1}^{\infty} \frac{nq^{2n}}{1-q^{2n}} \right) \\
&\equiv  \frac{(-q^2;q^2)_{\infty}}{(q^2;q^2)_{\infty}} \left( \sum_{n=1}^{\infty} \frac{(2n-1)q^{2n-1}}{1-q^{2n-1}} - \sum_{n=1}^{\infty} \frac{2nq^{2n}}{1-q^{2n}} \right) \pmod{3}\\
&=\frac{(-q^2;q^2)_{\infty}}{(q^2;q^2)_{\infty}}\sum_{n=1}^{\infty} \frac{(-1)^{n-1} n q^{n}}{1-q^{n}}\\
&\equiv \frac{(-q^2;q^2)_{\infty}}{(q^2;q^2)_{\infty}}\sum_{n=1}^{\infty} \frac{\chi(n) q^{n}}{1-q^{n}} \pmod{3},
\end{align*}
where $\chi(n)=1$ if $n\equiv 1$ or $2 \pmod{6}$, is $-1$ if $n\equiv 4$ or $5 \pmod{6}$, and is $0$ if $n\equiv 0 \pmod{3}$.
Thus,
\begin{align*}
M_1(q)& \equiv  \frac{(-q^2;q^2)_{\infty}}{(q^2;q^2)_{\infty}}\sum_{n=1}^{\infty}  E_{1,2} (n;6) q^n \pmod{3},
\end{align*}
where 
\begin{align*}
E_{1,2}(n;6)=\sum_{\substack{d\mid n\\ d\equiv 1,2\hspace{1mm}(\text{mod}\hspace{1mm}6)}} 1 \quad -\sum_{\substack{d\mid n\\ d\equiv -1,-2\hspace{1mm}(\text{mod}\hspace{1mm}6)}} 1.
\end{align*}
By \eqref{phi}, 
we see that
\begin{align*}
\phi(-q)&=\sum_{n=-\infty}^{\infty} (-1)^n q^{9n^2} -2q\sum_{n=-\infty}^{\infty} (-1)^n q^{9n^2+6n}\\
&=\phi(-q^9) -2q W(q^3). 
\end{align*}
Hence, by \cite[p. 80, eq. (32.39)]{fine},
{\allowdisplaybreaks
\begin{align*}
M_1(q) &\equiv \frac{1}{\phi(-q^2)} \sum_{n=1}^{\infty} E_{1,2} (n;6)q^n \pmod{3}\\
&\equiv \frac{1}{\phi(-q^2)} \left(1- \frac{(q^2;q^2)_{\infty} (q^3;q^3)_{\infty}^6}{(q;q)_{\infty}^2(q^6;q^6)_{\infty}^3} \right)  \pmod{3}\\
&=\frac{1}{\phi(-q^2)} \left(1- \frac{\phi(-q^3)^3}{\phi(-q)}\right)\\
&\equiv \frac{(\phi(-q)-\phi(-q^9))}{\phi(-q^2)\phi(-q)} \pmod{3}\\
&=\frac{-2qW(q^3) (-q;q)_{\infty} (-q^2;q^2)_{\infty}}{(q^2;q^2)_{\infty} (q;q)_{\infty}}\\
&=\frac{-2q W(q^3) (q;q)_{\infty}(-q^2;q^2)_{\infty}}{(q;q)_{\infty}^3}\\
&\equiv \frac{q W(q^3)}{(q^3;q^3)_{\infty}} (q^4;q^4)_{\infty} (q;q^2)_{\infty} \pmod{3}\\
&=\frac{q W(q^3)}{(q^3;q^3)_{\infty}} \sum_{n=-\infty}^{\infty} (-1)^n q^{2n^2-n}.
\end{align*}}
Now $2n^2-n$ is only congruent to $0$ or $1$ modulo $3$. Hence this last expression has no non-zero coefficients for terms where $q$ is a power of $3$. Hence $3\mid d_{3n}$. 

\subsection{Congruence modulo $6$}\label{cong6}
The congruence in \eqref{eq6n5} can be reduced to \eqref{equ_3n2}. 
Using \eqref{2mod}, we see that
\begin{align}\label{spt2cong}
\sum_{n=1}^{\infty}\overline{\text{spt}}_{\omega}(n)q^n&=\sum_{n=1}^{\infty}\frac{q^n(-q^{n+1};q)_n(-q^{2n+2};q^2)_{\infty}}{(1-q^n)^2(q^{n+1};q)_n(q^{2n+2};q^2)_{\infty}}\nonumber\\
&\equiv\sum_{n=1}^{\infty}\frac{q^n}{(1-q^n)^2}\pmod{2}\nonumber\\
&=\sum_{n=1}^{\infty}\sigma(n)q^n,
\end{align}
where $\sigma(n)$ denotes the sum of all positive divisors of $n$.

Now any number of the form $6n+5$ has all its prime divisors odd. The sum of the divisors of an odd prime raised to an odd power is even. For any number congruent to $5$ mod $6$ there must be at least one odd prime congruent to $5$ mod $6$ raised to an odd power in its prime factorization (otherwise the number would be congruent to $1$ mod $6$). Since $\sigma(n)$ is multiplicative, it follows that $\sigma(6n+5)$ is even. Hence the coefficients of $q^{6n+5}$ in both series are all even.

\subsection{Congruence modulo $5$}\label{cong5}
The congruence \eqref{equ_10n6} is proved here. Since $\overline{\textup{spt}}_2(5n+3)\equiv 0 \pmod{5}$, using \eqref{extra}, it suffices to show that $c_{10n+6} \equiv 0 \pmod{5}$, where $c_n$ is defined in \eqref{s22}. 
By \eqref{s22}, 
\begin{align*}
S(q)
&=\frac{(-q^2;q^2)_{\infty}}{(q^2;q^2)_{\infty}}\left(\sum_{n=1}^{\infty} \frac{(2n-1)q^{2n-1}}{1-q^{2n-1}} +2\sum_{n=1}^{\infty} \frac{(2n-1) q^{4n-2}}{1-q^{4n-2}}\right)\notag\\
&=\frac{(-q^2;q^2)_{\infty}}{(q^2;q^2)_{\infty}}\left(\sum_{n=1}^{\infty} \frac{(2n-1)q^{2n-1}}{1-q^{4n-2}} +\sum_{n=1}^{\infty} \frac{(2n-1)q^{4n-2}}{1-q^{4n-2}} +2\sum_{n=1}^{\infty} \frac{(2n-1) q^{4n-2}}{1-q^{4n-2}}\right)\notag\\
&= \frac{(-q^2;q^2)_{\infty}}{(q^2;q^2)_{\infty}}\left(\sum_{n=1}^{\infty} \frac{(2n-1)q^{2n-1}}{1-q^{4n-2}} +3 \sum_{n=1}^{\infty} \frac{(2n-1)q^{4n-2}}{1-q^{4n-2}} \right)\notag\\
&=: q E_1(q^2)+3E_2(q^2). 
\end{align*} 
Thus, it suffices to show that the coefficient of $q^{5n+3}$ in $E_2(q)$ is congruent to $0$ mod $5$, which follows from the following lemma. 

\newcommand{\QP}[1] {(#1)_{\infty}}

\begin{lemma}
Let
\begin{equation*}
 r(q) = q^{1/5} \frac{\QP{q,q^4;q^5}}{\QP{q^2,q^3;q^5}}\text{,}
\end{equation*}
and let $E_2(q)$ be defined as above. Then,
\begin{gather*}
E_{2}(q^{1/5}) \equiv \frac{q \QP{q;q}^2 \QP{q^{10};q^{10}}}{\QP{q^2,q^3;q^5}^5 \QP{q^5;q^5}^2} 
\left( \frac{r(q^2)}{r(q)^2} + \frac{1}{r(q)^2r(q^2)} + \frac{3}{r(q)^3} + \frac{r(q^2)}{r(q)^3} \right) \pmod{5}\text{.}
\end{gather*} 
\end{lemma}

\begin{proof}
 As in \cite{Tim}, set
\begin{equation*}
 A(q)= q^{1/5} \frac{\QP{q,q^4,q^5;q^5}}{\QP{q;q}^{3/5}}\text{,} \quad B(q)= \frac{\QP{q^2,q^3,q^5;q^5}}{\QP{q;q}^{3/5}}\text{.}
\end{equation*}
Although $A(q)^{\pm1}$ and $B(q)^{\pm1}$ do not have integer coefficients, all of the series $A(q)^{\pm5}$, $B(q)^{\pm5}$, and $r(q)^{\pm1}$ have integer coefficients. We will use the following properties:
\begin{align}
r(q) &= \frac{A(q)}{B(q)} \label{equ_rAB}\text{,}\\
A(q) B(q) &= q^{1/5} \frac{\QP{q^5;q^5}}{\QP{q;q}^{1/5}} \label{equ_ABqp}\text{,}\\
 2A(q)^5+B(q)^5 &\equiv 1 \pmod{5} \label{equ_A5B5}\text{,}\\
A(q^{1/5})^5 &\equiv \frac{r(q)}{1+2r(q)} \pmod{5} \label{equ_dis5}\text{.}
\end{align}
Identities \eqref{equ_rAB} and \eqref{equ_ABqp} follow directly from the definitions of $A(q)$ and $B(q)$. By multiplying through by $\QP{q;q}^3$ and using \eqref{equ_ABqp}, we see that \eqref{equ_A5B5} has the equivalent formulation
\begin{equation*}
 \QP{q^2,q^3,q^5;q^5}^5+2q \QP{q,q^4,q^5;q^5}^5 \equiv \QP{q;q}^3 \pmod{5}\text{.}
\end{equation*}
After applying Jacobi's triple product identity \cite[p.~21, Theorem 2.8]{gea1998} to each of the products on the left hand side, applying the fact that the characteristic is $5$, and using Jacobi's identity \cite[p.~176]{gea1998} on the right side, we must show that
\begin{equation*}
\sum_{n=-\infty}^{\infty} (-1)^n q^{\frac{5n}{2} (5 n-1)} +2q \sum_{n=-\infty}^{\infty} (-1)^n q^{\frac{5n}{2} (5 n-3)} \equiv \sum_{n=0}^{\infty} (2n+1)(-1)^{n} q^{\frac{n}{2} (n+1)} \pmod{5}\text{,}
\end{equation*}
which is easily seen to be true by breaking $n$ into residue classes modulo $10$ on the right hand side. Next, from \cite[Theorem 3.3]{Tim}, we have
\begin{equation*}
A(q^{1/5})^5 = A(q)^5-3 A(q)^4 B(q)+4 A(q)^3 B(q)^2-2 A(q)^2 B(q)^3+A(q) B(q)^4\text{,}
\end{equation*}
so, by \eqref{equ_rAB},
\begin{align*}
 A(q^{1/5})^5 &= A(q)^5-3 A(q)^4 B(q)+4 A(q)^3 B(q)^2-2 A(q)^2 B(q)^3+A(q) B(q)^4\\
&=B(q)^5 r(q) (1-2 r(q)+4 r(q)^2-3 r(q)^3 + r(q)^4)\\
& \equiv B(q)^5 r(q)(1+2r(q))^4 \pmod{5}\\
& \equiv \frac{r(q) (1+2r(q))^4}{1+2r(q)^5} \pmod{5}\\
& \equiv \frac{r(q)}{1+2r(q)} \pmod{5}\text{,}
\end{align*}
where \eqref{equ_A5B5} was used to obtain the penultimate equality. Thus, \eqref{equ_dis5} is clear. We will also require the two identities
\begin{align}
\QP{q^{1/5};q^{1/5}} &= q^{1/5} \QP{q^5;q^5} \left( \frac{1}{r(q)}-r(q)-1 \right) \label{equ_recip}\text{,}\\
 A(q)^5 &= \sum_{\substack{n=1\\ 5 \nmid n}}^{\infty} \frac{q^n}{1-q^n} \cdot \small \begin{cases}
                                                      1, &  n \equiv 1 \pmod{5}\\
                                                      -3, &  n \equiv 2 \pmod{5}\\
                                                      3, &  n \equiv 3 \pmod{5}\\
                                                      -1, &  n \equiv 4 \pmod{5}
                                                     \end{cases} \normalsize
\label{equ_eisen1}\text{.}
\end{align}
The first identity can be found in \cite[p. 270]{Part3} and the second is the first equality in \cite[Lemma 2.4]{Tim}. Let us save space by writing $r=r(q)$ and $R=r(q^2)$. Now, by \eqref{equ_eisen1},
\begin{align*}
E_2(q^{1/5})= \frac{\QP{-q^{1/5};q^{1/5}}}{\QP{q^{1/5};q^{1/5}}} & \sum_{n=1}^{\infty} \frac{(2n-1)q^{(2n-1)/5}}{1-q^{(2n-1)/5}} = \frac{\QP{q^{2/5};q^{2/5}}}{\QP{q^{1/5};q^{1/5}}^2} \sum_{n=1}^{\infty} \left(\frac{nq^{n/5}}{1-q^{n/5}} - \frac{2nq^{2n/5}}{1-q^{2n/5}} \right)\\
&\equiv \frac{\QP{q^{2/5};q^{2/5}}\QP{q^{1/5};q^{1/5}}^3}{\QP{q;q}} \sum_{n=1}^{\infty} \left(\frac{nq^{n/5}}{1-q^{n/5}} - \frac{2nq^{2n/5}}{1-q^{2n/5}} \right) \pmod{5} \\
&\equiv \frac{\QP{q^{2/5};q^{2/5}}\QP{q^{1/5};q^{1/5}}^3}{\QP{q;q}} \left(A(q^{1/5})^5 - 2 A(q^{2/5})^5 \right) \pmod{5} \text{.}
\end{align*}
By the dissection formulas \eqref{equ_recip} and \eqref{equ_dis5},
\begin{align*}
 E_2(q^{1/5}) &\equiv \frac{q \QP{q^{5};q^{5}}^3 \QP{q^{10};q^{10}}}{\QP{q;q}} \left(\frac{1}{R}-R-1\right)\left(\frac{1}{r}-r-1\right)^3 \left(\frac{r}{1+2r}-\frac{2 R}{1+2R}\right) \pmod{5} \\
&\equiv \frac{q \QP{q^{5};q^{5}}^3 \QP{q^{10};q^{10}}}{\QP{q;q}} \frac{(1+2 R)^2}{R} \frac{(1+2 r)^6}{r^3} \left(\frac{r}{1+2r}-\frac{2 R}{1+2R}\right) \pmod{5} \\
&\equiv \frac{q \QP{q^{5};q^{5}}^3 \QP{q^{10};q^{10}}(1+2 r^5)}{\QP{q;q}} \frac{(1+2 R)^2}{R} \frac{1+2 r}{r^3} \left(\frac{r}{1+2r}-\frac{2 R}{1+2R}\right) \pmod{5}\\
&= \frac{q \QP{q^{5};q^{5}}^3 \QP{q^{10};q^{10}}(1+2 r^5)}{\QP{q;q}} \left(-\frac{4 R}{r^2}+\frac{1}{r^2 R}-\frac{2}{r^3}-\frac{4 R}{r^3}\right) \\
&\equiv \frac{q \QP{q^{5};q^{5}}^3 \QP{q^{10};q^{10}}}{\QP{q;q}B(q)^5} \left(\frac{R}{r^2}+\frac{1}{r^2 R}+\frac{3}{r^3}+\frac{R}{r^3}\right) \pmod{5}\text{.}
\end{align*}
This is the result as stated.
\end{proof}

\noindent {\bf Remark.}
In the proof of Theorem~6.4 in \cite{ady1}, it is written that 
\begin{align*}
\sum_{n=0}^{\infty} \frac{(2n+1)q^{2n+1}}{1-q^{2n+1}}=q\frac{(q^4;q^4)_{\infty}^8}{(q^2;q^2)_{\infty}^4},
\end{align*}
which is not correct. What the authors meant is
\begin{align*}
\sum_{n=0}^{\infty} \frac{(2n+1)q^{2n+1}}{1-q^{4n+2}}=q\frac{(q^4;q^4)_{\infty}^8}{(q^2;q^2)_{\infty}^4}.
\end{align*}
The proof can be easily fixed as the congruence concerns the odd power terms only. However, using the functions $A(q)$, $B(q)$ and $r(q)$ from the proof of the lemma above, we can also correct the proof of Theorem~6.4 in \cite{ady1}. With the parameter $k= r(q) r(q^2)^2$, we have the parameterizations \cite[Entry 24]{contfrac}
\begin{equation*}
 r(q)^5 = k \left( \frac{1-k}{1+k} \right)^2 \text{,} \quad r(q^2)^5 = k^2 \left( \frac{1-k}{1+k} \right)\text{,}
\end{equation*}
and the dissection of the relevant $q$-series is found to be
\begin{equation*}
\frac{1}{\QP{q^{1/5};q^{1/5}}} \sum_{n=1}^{\infty} {\frac{n q^{n/5}}{1-q^{n/5}}} \equiv \frac{r(q^2) B(q)^5 B(q^2)^5}{q^{2/5} \QP{q^{10};q^{10}}} \frac{(2+k)^3}{(1+k)^2} \left(4k+2r(q)(1+k) +r(q^2)\right) \pmod5. 
\end{equation*}

\section{Congruences involving $\overline{\textup{spt}}(n)$ and $\overline{\textup{spt}}_{\omega}(n)$}\label{congsptsptw}
This section is devoted to proving Theorems~\ref{thm2} and \ref{thm3}. We first need the following lemmas. 
\begin{lemma}
\label{lemma_lerch2eisen}
We have
\begin{align}
\sum_{n=1}^{\infty} \frac{q^{n(n+1)/2}}{1-q^n}&=\sum_{n=1}^{\infty} \frac{q^{2n-1}}{1-q^{2n-1}}.  \label{sq1}
\end{align}
\end{lemma}
\proof
This is proved in \cite[p.~28]{macmahon}.
\endproof

\begin{lemma}
We have
\begin{align}
\sum_{n=1}^{\infty} \frac{q^{2n-1}}{1-q^{2n-1}} &\equiv \sum_{n=1}^{\infty} (q^{n^2}+q^{2n^2}) \pmod{2}.  \label{sq2}
\end{align}
\end{lemma}

\proof
We have
\begin{align*}
\sum_{n=1}^{\infty} \frac{q^{2n-1}}{1-q^{2n-1}}&=\sum_{n=1}^{\infty} \frac{q^{n}}{1-q^{n}}-\sum_{n=1}^{\infty} \frac{q^{2n}}{1-q^{2n}}\nonumber\\
&\equiv \sum_{n=1}^{\infty}q^{n^2}+\sum_{n=1}^{\infty}q^{2n^2} \pmod{2},
\end{align*}
by \eqref{44}.
\endproof

\subsection{Proof of Theorem~\textup{\ref{thm2}}}

By \eqref{barspt} and \eqref{barsptomega}, we know that
\begin{align*}
\sum_{n=1}^{\infty} \overline{\textup{spt}}(n)q^n & \equiv \sum_{n=1}^{\infty} \overline{\textup{spt}}_{\omega} (n)q^n \pmod{2}\\& \equiv \sum_{n=1}^{\infty}\frac{nq^n}{1-q^n} \pmod{2}\\
&\equiv \sum_{n=1}^{\infty} \frac{q^{2n-1}}{1-q^{2n-1}} \pmod{2}.
\end{align*}
Therefore, the congruences in \eqref{equ_nm2} follow from \eqref{sq2}. 

\subsection{Proof of Theorem~\textup{\ref{thm3}}}

Let us introduce the series

\begin{equation*}
 T(q) = \sum_{n=1}^{\infty} q^{n^2}\text{.}
\end{equation*}
Several identities satisfied by this series are
\begin{align}
1+2T(-q) &= \frac{(q;q)_{\infty}}{(-q;q)_{\infty}} \label{equ_phi_prod}\text{,}\\
T(q)+T(q)^2 &= \sum_{\substack{n=0}}^{\infty} \frac{(-1)^{n}q^{2n+1}}{1-q^{2n+1}} \label{equ_two_squares}\text{,}
\end{align}
where the first identity is a restatement of \eqref{phi} and the second is \cite[p. 59, Equation (26.63)]{fine}. By applying \eqref{barspt}, 
\eqref{sq1}, \eqref{equ_phi_prod}, and \eqref{equ_two_squares}  in the same sequence, we find that
\begin{align*}
\sum_{n=1}^{\infty} \overline{\textup{spt}}(n)q^n &= \frac{(-q;q)_{\infty}}{(q;q)_{\infty}} \left( \sum_{n=1}^{\infty} \frac{nq^n}{1-q^n} + 2\sum_{n=1}^{\infty} \frac{(-1)^n q^{n^2+n}}{(1-q^n)^2} \right)\\
& \equiv \frac{(-q;q)_{\infty}}{(q;q)_{\infty}} \left( \sum_{\substack{n=0}}^{\infty} \frac{(-1)^{n}q^{2n+1}}{1-q^{2n+1}} + 2\sum_{\substack{n=0}}^{\infty} \frac{q^{4n+2}}{1-q^{4n+2}}  +2 \sum_{n=1}^{\infty} \frac{q^{n (n+1)}}{1-q^{2n}} \right) \pmod 4\\
& \equiv \frac{(-q;q)_{\infty}}{(q;q)_{\infty}} \sum_{\substack{n=0}}^{\infty} \frac{(-1)^{n}q^{2n+1}}{1-q^{2n+1}} \pmod 4\\
& \equiv \frac{1}{1+2T(-q)} (T(q)+T(q)^2) \pmod 4\\
& \equiv (1+2T(q))(T(q)+T(q)^2) \pmod 4\\
&\equiv  2T(q)^3+3T(q)^2+T(q) \pmod 4\text{.}
\end{align*}
Next, set
\begin{equation*}
\begin{alignedat}{5}
 t&=\sum_{n=1}^{\infty} q^{(7n+0)^2/7} =  T(q^7)\text{,}& \quad b&= \sum_{n=-\infty}^{\infty} q^{(7n+2)^2/7}\text{,}\\
a&= \sum_{n=-\infty}^{\infty} q^{(7n+1)^2/7}\text{,}& c&= \sum_{n=-\infty}^{\infty} q^{(7n+3)^2/7}\text{,}
\end{alignedat}
\end{equation*}
and note that $T(q^{1/7}) = t+a+b+c$. Therefore,
\begin{align*}
 \sum_{n=1}^{\infty} \overline{\textup{spt}}(n)q^{n/7} &\equiv  2T(q^{1/7})^3+3T(q^{1/7})^2+T(q^{1/7}) \pmod 4\\
&=12 a b c+2 t^3+3 t^2+t\\
&\quad+6 a t^2+6 a t+a+6 b^2 t+3 b^2+6 b c^2\\
&\quad+6 a^2 t+3 a^2+6 a b^2+6 c t^2+6 c t+c\\
&\quad+2 a^3+12 a c t+6 a c+6 b^2 c\\
&\quad+6 a^2 c+6 b t^2+6 b t+b+6 c^2 t+3 c^2\\
&\quad+12 a b t+6 a b+6 a c^2+2 b^3\\
&\quad+6a^2 b+12 b c t+6 b c+2 c^3\text{,}
\end{align*}
where the terms in the expansion have been collected according to the residue classes modulo $1$ of the exponents on $q$. Taking the terms involving only integral powers of $q$ gives
\begin{align*}
 \sum_{n=1}^{\infty} \overline{\textup{spt}}(7n)q^{n} &\equiv 12 a b c+2 t^3+3 t^2+t \pmod{4}\\
&\equiv 2 t^3+3 t^2+t \pmod{4}\\
&\equiv \sum_{n=1}^{\infty} \overline{\textup{spt}}(n)q^{7n}\pmod{4}\text{.}
\end{align*}
This proves the congruence \eqref{equ_7nm4}. Similarly, by \eqref{barsptomega}, 
{\allowdisplaybreaks\begin{align*}
\sum_{n=1}^{\infty} \overline{\textup{spt}}_{\omega}(n)q^n &= \frac{(-q^2;q^2)_{\infty}}{(q^2;q^2)_{\infty}} \left( \sum_{n=1}^{\infty} \frac{nq^n}{1-q^n} + 2\sum_{n=1}^{\infty} \frac{(-1)^n q^{2n(n+1)}}{(1-q^{2n})^2} \right)\\
& \equiv \frac{(-q^2;q^2)_{\infty}}{(q^2;q^2)_{\infty}} \left( \sum_{\substack{n=0}}^{\infty} \frac{(-1)^{n}q^{2n+1}}{1-q^{2n+1}} + 2\sum_{\substack{n=0}}^{\infty} \frac{(-1)^n q^{4n+2}}{1-q^{4n+2}}  +2 \left( \sum_{n=1}^{\infty} \frac{q^{n (n+1)/2}}{1-q^{n}} \right)^4 \right) \pmod 4\\
& = \frac{1}{1+2T(-q^2)} \left( T(q)+T(q)^2+2T(q^2)+2T(q^2)^2 + 2( T(q)+T(q)^2 )^4 \right) \pmod 4\\
& \equiv (1+2T(q)^2)(T(q)+ 3T(q)^2+2T(q)^8 ) \pmod 4\\
& \equiv 2 T(q)^8+ 2T(q)^4+ 2T(q)^3+ 3 T(q)^2+T(q) \pmod 4\text{.}
\end{align*}}
This polynomial in $T(q)$ may be dissected in the same fashion, and we find that
\begin{align*}
 \sum_{n=1}^{\infty} \overline{\textup{spt}}_{\omega}(7n)q^{n} &\equiv 2 t^8+ 2t^4+2t^3+ 3 t^2+t \pmod 4\\
&\equiv \sum_{n=1}^{\infty} \overline{\textup{spt}}_{\omega}(n)q^{7n} \pmod 4\text{,}
\end{align*}
which proves the congruence \eqref{equ_w7nm4}.

\section{Concluding remarks and some open problems}\label{conclude}
In concluding, we remark that the study involving the partition function $\overline{p}_{\omega}(n)$ and its associated smallest parts function $\overline{\textup{spt}}_{\omega}(n)$ is more difficult than the one involving $p_{\omega}(n)$ and $\textup{spt}_{\omega}(n)$. Nonetheless, these functions are as interesting as their aforementioned counterparts. This certainly merits further study of these functions. In particular, we give two open problems below.\\

\textbf{Problem 1.} Is the generating function of $\overline{p}_{\omega}(n)$ representable in terms of some more important functions, say for example, mock theta functions or more generally, mock modular forms or mixed mock modular forms? Is it possible to further simplify the series involving the little $q$-Jacobi polynomials in Theorem \ref{mt1}?\\

We were able to obtain another proof of the mod $4$ congruences in Theorem \ref{pw4} assuming the conjecture that the function $Y(q)$ defined by 
\begin{equation}\label{yq}
Y(q):=\sum_{n,m\geq 1}\frac{(-1)^mq^{2nm+m}}{(1+q^n)(1-q^{2m-1})}
\end{equation}
is an odd function of $q$. Some coefficients in the expansion of $Y(q)$ are
\begin{equation*}
Y(q)=-q^3-2q^5-3q^7-5q^9-4q^{11}-7q^{13}-9q^{15}-\cdots-53q^{91}-62q^{93}-38q^{95}-55q^{97}-\cdots.
\end{equation*}
Unfortunately we are unable to prove that indeed it is an odd function, hence we state it below as another open problem.

\textbf{Problem 2.} Prove that the function $Y(q)$ defined in \eqref{yq} is an odd function of $q$.

\begin{center}
\textbf{Acknowledgements}
\end{center}
The fourth author was partially supported by a grant ($\#280903$) from the Simons Foundation.


\begin{thebibliography}{00}
\bibitem{agarwal}
R.P.~Agarwal, \emph{A family of basic hypergeometric and combinatorial identities and certain summation formulae}, Indian J. Pure Appl. Math.~\textbf{12} (1981), 728--737.

\bibitem{agar1}
R.P.~Agarwal, \emph{On the paper ``A `Lost' Notebook of Ramanujan''}, Adv. Math.~\textbf{53} (1984), 291--300.

\bibitem{alladi}
K.~Alladi, \emph{Variants of classical q-hypergeometric identities and partition implications}, Ramanujan J.~\textbf{31} Issues 1-2 (2013), 213--238.

\bibitem{gea55} 
G.E.~Andrews, \emph{On the $q$-analog of Kummer's theorem and applications}, Duke Math. J.~\textbf{40} (1973), 525--528.

\bibitem{gea1998}
G.E.~Andrews, The theory of partitions, Addison-Wesley Pub. Co., NY, 300 pp. (1976). Reissued, Cambridge University Press, New York, 1998.

\bibitem{gea90}
G.E.~Andrews, \emph{Ramanujan's ``Lost" notebook: \textup{I}. Partial theta functions}, Adv. Math.~\textbf{41} (1981), 137--172.

\bibitem{andrews1984}
G.E.~Andrews, \emph{Ramanujan's ``Lost'' Notebook \textup{IV}. Stacks and Alternating Parity in Partitions}, Adv. Math.~\textbf{53} (1984), 55--74.

\bibitem{aa}
G.E.~Andrews and R.A.~Askey, \emph{Enumeration of partitions: The role of Eulerian series and $q$-orthogonal polynomials}, Higher Combinatorics, M. Aigner, ed., Reidell Publ. Co., Dordrecht, Holland, pp. 3-26 (1977).

\bibitem{contfrac} G.E. Andrews, B.C. Berndt, L. Jacobsen, and R.L. Lamphere, The continued fractions found in the unorganized portions of Ramanujan's notebooks, Memoir No. 477, Amer. Math. Soc. 99 (1992).

\bibitem{ady1}
G.E.~Andrews, A.~Dixit and A.J.~Yee, \emph{Partitions associated with the Ramanujan/Watson mock theta functions $\omega(q), \nu(q)$ and $\phi(q)$}, \emph{Research in Number Theory},~\textbf{1} Issue 1 (2015), 1--25.

\bibitem{bailey}
W.N.~Bailey, \emph{Identities of the Rogers-Ramanujan type}, Proc. London Math. Soc.~\textbf{50} No. 2 (1949), 1--10.

\bibitem{bcgkmw}
J.~Berg, A.~Castillo, R.~Grizzard, V.~Kala, R.~Moy and C.~Wang, \emph{Congruences for Ramanujan's $f$ and $\omega$ functions via generalized Borcherds products}, Ramanujan J.~\textbf{35} (2014), 327--338.

\bibitem{Part3}  B. C. Berndt, Ramanujan's Notebooks, Part III, Springer-Verlag, New York, 1991.

\bibitem{bringmannlovejoyosburn1}
K. Bringmann, J. Lovejoy, and R. Osburn, \emph{Rank and crank moments for overpartitions}, J. Number Theory, \textbf{129} (2009), 1758--1772.

\bibitem{bringmannlovejoyosburn2}
K. Bringmann, J. Lovejoy, and R. Osburn, \emph{Automorphic properties of generating functions for generalized rank moments and Durfee symbols}, IMRN \textbf{2010} No. 2, 238--260.

\bibitem{jbko}
J.~Bruinier and K.~Ono, \emph{Identities and congruences for Ramanujan's $\omega(q)$}, Ramanujan J.~\textbf{23} (2010), 151--157.

\bibitem{corlove}
S.~Corteel and J.~Lovejoy, \emph{Overpartitions}, Trans. Amer. Math. Soc.~\textbf{356} No. 4 (2004), 1623--1635.

\bibitem{fine}
N.J. Fine, {Basic hypergeometric series and applications}, Amer. Math. Soc., Providence, 1988.

\bibitem{garthwaite}
S.A.~Garthwaite, \emph{The coefficients of the $\omega(q)$ mock theta function}, Int.~J.~Number Theory~\textbf{4} No. 6 (2008), 1027--1042.

\bibitem{jennings-shaffer}
F.G. Garvan and C. Jennings-Shaffer, \emph{The spt-crank for overpartitions}, Acta Arithmetica, \textbf{166} (2014), 141--188.

\bibitem{gasper} 
G.~Gasper and M.~Rahman, {Basic Hypergeometric Series}, 2nd ed., Cambridge University Press, Cambridge, 2004.

\bibitem{Tim} T.~Huber, \emph{A theory of theta functions to the quintic base}, J.  Number Theory, \textbf{134} (2014), 49--92.

\bibitem{ismail}
M.E.H.~Ismail, {Classical and Quantum Orthogonal Polynomials in One Variable}, Cambridge University Press, Cambridge, 2009.

\bibitem{iz}
M.E.H.~Ismail and R.~Zhang, \emph{$q$-Bessel functions and Rogers-Ramanujan type identities}, arXiv:1508.06861v1.

\bibitem{macmahon}
P.A.~Macmahon, \emph{Combinatory Analysis}, Vol.~II, Cambridge University Press, Cambridge, 1915-1916, reissued, Chelsea, 1960.

\bibitem{robbins}
N.~Robbins, \emph{On partition functions and divisor sums}, J. Integer Sequences~\textbf{5} (2002), Article 02.1.4.

\bibitem{waldherr}
M.~Waldherr, \emph{On certain explicit congruences for mock theta functions}, Proc.~Amer.~Math.~Soc.~\textbf{139} No. 3 (2011), 865--879.
\end{thebibliography}
\end{document}